\input amstex
\documentstyle{amsppt}
\magnification=1200
\hcorrection{0.3in}
\vcorrection{-0.2in}


\define\op#1{\operatorname{#1}}

\NoRunningHeads
\NoBlackBoxes
\topmatter
\title A Geometric Approach to the Modified Milnor Problem
\endtitle
\author Lina Chen, \footnote{Supported partially by China Postdoctoral Science Foundation 2017M620143.
\hfill{$\,$}} Xiaochun Rong \footnote{Supported by a research fund from Capital Normal University.\hfill{$\,$}}
and Shicheng Xu\footnote{Supported partially by NSFC
Grant 11401398 and by Youth Innovative Research Team of Capital Normal University.\hfill{$\,$}}
\endauthor

\address Mathematics Department, East China Normal University, Shanghai, P.R.C.
\endaddress
\email chenlina\_mail\@163.com
\endemail
\address Mathematics Department, Rutgers University
New Brunswick, NJ 08903 USA
\endaddress
\email rong\@math.rutgers.edu
\endemail
\address Mathematics Department, Capital Normal University, Beijing,
P.R.C.
\endaddress
\email  shichxu\@foxmail.com
\endemail

\abstract 
The Milnor Problem (modified) in the theory of group growth asks whether any finite presented group of vanishing algebraic entropy has at most polynomial growth. We show that a positive answer to the Milnor Problem (modified) is equivalent to the Nilpotency Conjecture in Riemannian geometry: given $n, d>0$, there exists a constant $\epsilon(n,d)>0$ such that if a compact Riemannian $n$-manifold $M$ satisfies that Ricci curvature $\op{Ric}_M\ge -(n-1)$, diameter $d\ge \op{diam}(M)$ and volume entropy $h(M)<\epsilon(n,d)$,  then the fundamental group $\pi_1(M)$ is virtually nilpotent.  We will verify the Nilpotency Conjecture in some cases, and we will verify the vanishing gap phenomena for more cases i.e., if $h(M)<\epsilon(n,d)$, then $h(M)=0$.
\endabstract
\endtopmatter
\document

\head 0. Introduction
\endhead

\vskip4mm

Given a finitely generated group $\Gamma$ and a symmetric finite generating set $S$ i.e.
$\gamma\in S$ if and only if $\gamma^{-1}\in S$, for any $\gamma\in \Gamma$, its
word length is defined by
$$|\gamma|_S=\min \{k,\,\, \gamma=\gamma_{i_1}\cdots \gamma_{i_k},\,\, \gamma_{i_j}\in S\}.$$
Let $|\Gamma(R,S)|=\#\{\gamma\in \Gamma, |\gamma|_S\le R\}$. The algebraic entropy of $(\Gamma,S)$ is an invariant measuring the exponential growth rate of $\Gamma$, defined by
$$h_w(\Gamma,S)=\lim_{R\to \infty}\frac {\ln |\Gamma(R,S)|}R,$$
which always exists (see Section 1).
We say that $\Gamma$ has at most polynomial growth,
if $|\Gamma(R,S)|\le R^m$ for all large $R$ and some constant $m$; the property of at most polynomial growth is independent of $S$. Clearly, any  finitely generated group of at most polynomial growth has a vanishing algebraic
entropy, $h_w(\Gamma,S)=0$.

In the theory of group growth, a well-known Milnor's problem (\cite{Mi3}) asks whether vanishing of $h(\Gamma,S)$ implies that $\Gamma$ is at most polynomial growth. A negative answer was obtained with a finitely generated but not finitely presented group by Grigorchuk
\cite{Gri1-3}, where he modified Milnor's problem as:

\example{Milnor Problem 0.1} Let $\Gamma$ be a finitely presented group.
If $h_w(\Gamma,S)=0$, then $\Gamma$ has at most polynomial growth.
\endexample

It is conjectured in \cite{GP} that Milnor Problem 0.1 has a positive answer. For a brief history on
Milnor Problem 0.1, see \cite{Gri4}.

To a geometer, the significance of Milnor Problem 0.1 is that $\Gamma$  can be treated as the fundamental
group of a compact $n$-manifold with $n\ge 4$ (\cite{CZ}), a major subject in Riemannian geometry
with a rich history (cf. \cite{Pe}).

The purpose of this paper is to propose the following nilpotency conjecture on compact $n$-manifolds
with bounded Ricci curvature and diameter, and prove its equivalence to a positive
answer to Milnor Problem 0.1. It is our hope that this equivalence will
bring geometric tools into the study of Milnor Problem 0.1.

Let $M$ be a compact Riemannian $n$-manifold. The volume entropy of $M$, $h(M)$,
is an asymptotic geometric invariant that measures the exponential growth rate of the volume of metric balls in the Riemannian universal cover $\tilde M$. Precisely,
$$h(M)=\lim_{R\to \infty}\frac{\ln \op{vol}(B_R(\tilde p))}R,$$
where the limit always exists and is independent of $\tilde p\in \tilde M$ (\cite{Man}). Let $\Gamma=\pi_1(M)$.
It is not hard to show that (\cite{Mi1})
$$C^{-1}\cdot h_w(\Gamma,S)\le h(M)\le C\cdot h_w(\Gamma,S),$$
which implies that $h_w(\Gamma,S)=0$ if and only if $h(M)=0$, where constant $C$ may depends on $S$.

\example{Nilpotency Conjecture 0.2} Given $n, d>0$, there exists a constant $\epsilon(n,d)>0$, such that if a complete  $n$-manifold $M$ satisfies
$$\op{Ric}_M\ge -(n-1),\quad d\ge \op{diam}(M),\quad h(M)<\epsilon(n,d),$$
then $\pi_1(M)$ contains a nilpotent subgroup of finite index (i.e., $\pi_1(M)$ is virtually nilpotent).
\endexample

Here is a partial motivation for Nilpotency Conjecture 0.2: a result in Riemannian geometry, which was conjectured by Gromov and proved by Kapovitch-Wilking (\cite{KW}), says that there exist constants,  $\epsilon(n), c(n)>0$,
such that if a complete $n$-manifold $M$ satisfies $\op{Ric}_M\ge -(n-1)$ and $\op{diam}(M)\le \epsilon(n)$, then $\pi_1(M)$ has a nilpotent subgroup of index $\le c(n)$ (briefly,  $c(n)$-nilpotent).  A natural question is when replacing ``$\op{diam}(M)\le \epsilon(n)$'' with $\op{diam}(M)\le d$ for any $d>0$,
what additional geometric condition on $M$ is required to conclude a virtual nilpotency for $\pi_1(M)$?
Nilpotency Conjecture 0.2 provides with an answer.

The main result of this paper is

\proclaim{Theorem 0.3} Nilpotency Conjecture 0.2 is equivalent to a positive answer to Milnor Problem 0.1.
\endproclaim

Perhaps Theorem 0.3 is a little surprise, because Milnor Problem 0.1 and Nilpotency
Conjecture 0.2 arise from different fields and are apparently not related. By Theorem 0.3,
progresses made on either problem will shed a light on the other.

\remark{Remark \rm 0.4}  A further problem on Nilpotency Conjecture 0.2 is whether $\pi_1(M)$
can be strengthened to be $w(n,d)$-nilpotent. In the case that $\pi_1(M)$ is torsion free,  we will
strengthen Theorem 0.3 to that $\pi_1(M)$ is $w(n,d)$-nilpotent (see Theorem 7.3), provided a positive answer to Milnor Problem 0.1. We point out that our argument may not rule out virtually nilpotent torsion groups with arbitrarily large index, see Example 7.4.
\endremark

We now explain our approach to Theorem 0.3. First, for $n\ge 4$, any finitely presented group  can be realized as $\pi_1(M)$ for some compact $n$-manifold $M$ (Theorem 5.1.1 in \cite{CZ}), and
$h_w(\pi_1(M),S)=0$ iff $h(M)=0$. Because a virtually nilpotent group has at most  polynomial growth (see Lemma 1.4, cf. [Wo]), Nilpotency Conjecture 0.2 implies a positive answer to Milnor  Problem 0.1.
The main work in Theorem 0.3 is to show that $h(M)<<1$ implies that $\pi_1(M)$ is virtually nilpotent, assuming a positive answer to Milnor Problem 0.1.

By the Gromov's compactness, it is equivalent to show, for a convergent sequence of compact $n$-manifolds, $M_i@>GH>>X$, satisfying
$$\op{Ric}_{M_i}\ge -(n-1),\quad d\ge \op{diam}(M_i), \quad h(M_i)\to 0,\tag 0.5 $$
that $\Gamma_i=\pi_1(M_i)$ is virtually nilpotent for
$i$ large. According to \cite{Gr}, it reduces to show that $\Gamma_i$ has at most polynomial growth.

According to [FY], passing to a subsequence the sequence of Riemannian
universal cover equivariantly converges, $(\tilde M_i,\tilde p_i,\Gamma_i)@>GH>>(\tilde X,\tilde p,G)$, and there is $\epsilon>0$ such that the subgroup $\Gamma_{i,\epsilon}=\{\gamma_i\in \Gamma_i,\, d(\tilde q_i,\gamma_i(\tilde q_i))<\epsilon, \,\, \tilde q_i\in B_1(\tilde p_i)\}$ is normal and for $i$ large, $\Gamma_i/\Gamma_{i,\epsilon}\cong G/G_0$, where $G_0$ is the identity
component of $G$. Let's divide the rest into two situations: non-collapsing i.e,
$\op{vol}(M_i)\ge v>0$ or collapsing i.e., $\op{vol}(M_i)\to 0$ as $i\to \infty$.

In the non-collapsed situation, $G_0=e$, $\Gamma_{i,\epsilon}=e$ i.e., $\Gamma_i=G$ for $i$ large (see Lemma 1.10, cf. \cite{PR}). Consequently, $h_w(\Gamma_i)=0$ and thus $\Gamma_i$ has at most polynomial growth (by a positive answer to Milnor Problem 0.1).

Assuming $\op{vol}(M_i)\to 0$, we apply the stable structural result on fundamental groups in \cite{KW} (see Theorem 1.12) to $(\Gamma_i,\Gamma_{i,\epsilon})$, and conclude that $\Gamma_{i,\epsilon}$ contains a nilpotent subgroup, $N_i$, which is normal in $\Gamma_i$ such that $\Lambda_i=\Gamma_i/N_i$ has a bounded number of possible isomorphic classes, and
there is a descending sequence of a uniform length $k\le n$,
$$1\to N_i\to \Gamma_i @>\pi_i >> \Lambda_i\to 1,\quad N_i=N_{i1}\vartriangleright\cdots \vartriangleright N_{ik}\vartriangleright N_{ik+1}=\operatorname{Tor}(N_i),\tag 0.6$$
satisfies that $[N_i,N_{ih}]\subseteq N_{ih+1}$, $N_{ih}/N_{ih+1}$ is
a free abelian ($1\le h\le k$), and the induced maps
via conjugation, $\bar\rho_h: \Lambda_i\to \op{Aut}(N_{ih}/N_{ih+1})$,
admit only finite choices up to an automorphism, where $\op{Tor}(N_i)$
denotes the torsion subgroup of $N_i$.

In view of the above, passing to a subsequence we may assume that $\Lambda_i=\Lambda$ and $\bar \rho_h: \Lambda\to \op{Aut}(N_{ih}/N_{ih+1})$ are independent of large $i$. Unless specified elsewhere, we will always take a symmetric
generating set, $S_i=B_{i}\cup S_{0,i}$, such that $B_i\cap S_{0,i}=\emptyset$, $B_i$ generates $N_i$ and $\bar S_{0}=\pi_i(S_{0,i})$ (independent of $i$) generates $\Lambda_i$.

Our main effort is to establish the following algebraic result
which provides a sufficient condition for $\Gamma$ to be at most polynomial growth.

\proclaim{Theorem 0.7} A finitely generated group $\Gamma$ has at most polynomial growth, if $\Gamma$ has a normal nilpotent subgroup $N$ satisfying
(0.6) and the following two conditions:

\noindent \rm{(0.7.1)} $\Lambda$ is virtually nilpotent, and $\bar S_0$ is a graded symmetric generating set of $\Lambda$.

\noindent \rm{(0.7.2)} For $\gamma\in S_0$, if $\pi(\gamma)$ has an infinite order, then all eigenvalues of $\bar \rho_h(\pi(\gamma))$ have norm one, $1\le h\le k$.
\endproclaim

The proof of Theorem 0.7 is somewhat tedious and technical (see Section 3
for an outline), which is similar to the approach in \cite{Wo}.

To conclude Theorem 0.3 from Theorem 0.7, it remains to check (0.7.1) and (0.7.2) for $(\Gamma_i,N_i)$. Then (0.7.1) follows from $h_w(\Lambda,\bar S_0)\le h_w(\Gamma_i,S_i)\to 0$ and a positive answer of Milnor Problem 0.1. A verification of (0.7.2) is based on a lower bound estimate on $h_w(\Gamma)$, $\Gamma$ satisfies that $1\to \Bbb Z^k\to \Gamma\to \Bbb Z\to 1$ (\cite{Os1}), and the finiteness of $\bar \rho_h$ (see Lemma 2.1).

Nilpotency Conjecture 0.2 has been known in a few cases; a consequence of the splitting theorem of Cheeger-Gromoll (\cite{CG}) is that the fundamental group of any compact manifold of non-negative Ricci curvature is virtually abelian. Nilpotency Conjecture 0.2 holds in dimension $n\le3$; if $\Gamma$ is infinite, then $\Gamma$ is either of
polynomial growth or its uniform exponential growth rate $\inf_{S} \{h_w(\Gamma,S)\}$ is bounded from below by a universal positive constant \cite{Ce}
(see Lemma 1.2).

Using results from \cite{Ch} and \cite{Ti}, we will verify Nilpotency Conjecture 0.2 for more cases by showing that (0.7.1) holds. A discrete group $\Gamma$ is called amenable, if $\Gamma$ has a finite additive left-invariant probability
measure. An elementary amendable group is the smallest class of amendable groups containing finite groups,
abelian groups and preserved by operations: subgroups, quotient groups, groups extension and direct limit.

\proclaim{Theorem 0.8} Nilpotency Conjecture 0.2 holds in the following cases:

\noindent \rm{(0.8.1)} $\pi_1(M)$ is elementary amenable.

\noindent \rm{(0.8.2)} $\pi_1(M)$ is linear, i.e., a subgroup
of $\op{GL}(m,\Bbb R)$.
\endproclaim

The Nilpotency Conjecture 0.2 implies a vanishing gap phenomena for $h(M)$ that $h(M)<\epsilon(n,d)$ implies that $h(M)=0$. Contrasting to the gap vanishing
phenomena (Theorem 0.8), $\{h_w(\Gamma,S), \, \Gamma \text{ is elementary
amendable}\}$ has ``$0$'' as an accumulation point (\cite{Os2}).

The vanishing gap phenomena always holds for non-collapsed $M$, i.e.
$\op{vol}(M)\ge v>0$, with replacing $\epsilon(n,d)$ by $\epsilon(n,d,v)$ (Lemma 6.2).
We will verify the gap phenomena in more cases in which
the Nilpotency Conjecture 0.2 are open.

\proclaim{Theorem 0.9} Let $M$ be a compact $n$-manifold such that
$$\op{Ric}_M\ge -(n-1),\quad d\ge \op{diam}(M),$$
and $\Gamma$ contains a normal nilpotent subgroup $N$ satisfying (0.6).
Then $h(M)=0$ or $h(M)\ge \epsilon(n,d)$ under either of
the following conditions:

\noindent \rm{(0.9.1)} $\Gamma\cong N\rtimes_{\psi} \Lambda$.

\noindent \rm{(0.9.2)} $\Lambda$ satisfy a polynomial isoperimetric inequality (see Sec. 5 for the definition).
\endproclaim

In the proof of Theorem 0.9, a key is to establish a result similar to Theorem 0.7
(see Theorem 5.4); if replacing (0.7.1) by that $\op{Im}(\bar \rho_h)$ is virtually nilpotent, then $h_w(\Gamma,S)=h_w(\left <S_0\right>,S_0)$, or $h_w(\Gamma,S)=h_w(\Lambda,\bar S_0)$ under (0.9.1) or (0.9.2), where $\left<S_0\right>$ denotes the subgroup generated by elements in $S_0$.

\vskip4mm

The rest of the paper is organized as follows:

\vskip2mm

In Section 1, we will supply basic notions and properties required in the rest of the paper.

In Section 2, we will prove Theorem 0.3 and Theorem 0.8 by assuming theorem 0.7.

In Section 3, we will outline the proof of Theorem 0.7.

In Section 4, we will prove Theorem 0.7.

In Section 5, we will prove Theorem 0.9.

In Section 6, we will show that $h(M)$ small implies that $\pi_1(M)$ is
amendable, and we will verify a few cases of gap vanishing volume entropy.

In Section 7, we will show that in Theorem 0.3, when $\pi_1(M)$ is torsion free, then $\pi_1(M)$ is $w(n,d)$-nilpotent.

\vskip4mm {\bf Acknowledgement}. The authors would like to thank S. Honda
for asking the third author the question related to \cite{CRX} whether there
is a vanishing volume entropy gap phenomena, and for several useful
discussions with the third author.

\vskip10mm

\head 1. Preliminaries
\endhead

\vskip4mm

In this section, we will supply notions and properties that will
be used in the proof of Theorems 0.3, 0.7, 0.8 and 0.9.

\vskip4mm

\subhead a. Algebraic, geometric and volume entropies
\endsubhead

\vskip4mm

Let $M$ be a compact manifold, let  $\pi: (\tilde M,\tilde p)\to
(M,p)$ denote the Riemannian universal cover, let $\Gamma=\pi_1(M,p)$ be
the fundamental group, and let $S$ be a finite symmetric generating set
 for $\Gamma$. The definitions of $h_w(\Gamma,S)$ and $h(M)$
were given in the introduction without a justification on the existence of
limits.

For the convenience of readers, we will give a brief account on the existence of limits. For this purpose, we introduce the notion of geometric entropy of $\Gamma$:
$$h(\Gamma)=\lim_{R\to \infty}\frac{\ln |\Gamma(R)|}R,$$
where $\Gamma(R)=\{\gamma\in \Gamma, \,\, |\tilde p\gamma(\tilde p)|<R\}$.
The existence of the above limit is based on the property: for any $r, s>> d=\op{diam}(M)$,
$$|\Gamma(r+s)|\le |\Gamma(r)|\cdot |\Gamma(s+2d)|.$$
For $R>>r>>d$ and $(k-1)r\le R\le kr$, iterating the above inequality $k$ times yields
that $\ln |\Gamma(R)|\le k\cdot \ln |\Gamma(r+2d)|$, and thus
$$\limsup_{R\to \infty}\frac{\ln |\Gamma(R)|}R\le \frac {\ln |\Gamma(r+2d)|}r.$$
Consequently,
$$\limsup_{R\to \infty}\frac{\ln |\Gamma(R)|}R\le \liminf_{r\to \infty}\frac {\ln |\Gamma(r+2d)|}r=
\liminf_{r\to \infty}\frac{\ln |\Gamma(r+2d)|}{r+2d},$$
which proves the existence of $h(\Gamma)$.

On the other hands, for any $R>>d$ it is easy to check,
$$\frac{\op{vol}(B_{R-d}(\tilde p))}{\op{vol}(B_d(\tilde p))}\le |\Gamma(R)|\le \frac{\op{vol}(B_{R+d}(\tilde p))}{\op{vol}(M)},$$
which implies that $h(\Gamma)$ exists if and only if $h(M)$ exists. Summarizing the above, we have

\proclaim{Lemma 1.1 (\cite{Man})} Let $M$ be a compact manifold with fundamental group $\Gamma=\pi_1(M,p)$. Then
$$h(M)=\lim_{R\to\infty}\frac{\ln \op{vol}\,(B_R(\tilde p))}R,$$
exists and $h(M)=h(\Gamma)$.
\endproclaim

To see the existence of the limit for $h_w(\Gamma)$, one may start with the inequality,
$$|\Gamma(r+s,S)|\le |\Gamma(r,S)|\cdot |\Gamma(s,S)|,$$
and by a similar discussion, one gets the existence of the following limit,
$$h_w(\Gamma,S)=\lim_{R\to \infty}\frac{\ln |\Gamma(R,S)|}R.$$


A basic relation between $h_w(\Gamma,S)$ and $h(M)$ is that there is a constant
$C>0$ depending on $S$ and $M$ such that $C^{-1}\cdot h(M)\le h_w(\Gamma,S)
\le C\cdot h(M)$ (\cite{Mil}). Let $d=\op{diam}(M)$. Note that $\Gamma(2d+1)$
is symmetric and finite. For any $\gamma\in \Gamma$,
by dividng $\gamma$ into pieces of length $r<1$, $\gamma$ can be expressed
as a product of elements in $\Gamma(2d+1)$.

\proclaim{Lemma 1.2} Let $M$ be a compact manifold of diameter $d$.
For $p\in M$ and $S=\Gamma(2d+1)$,
$$h(M)\le h_w(\Gamma,S)\le (2d+1)h(M).$$
\endproclaim

\demo{Proof} By Lemma 1.1, it remains to show that $\Gamma(R)\subseteq \Gamma(R,S)\subseteq \Gamma((2d+1)R)$.
Note that the right hand inclusion is clear, due to the selection for $S$.

For any $\gamma\in \Gamma(R)$, we may assume the lifting, $\tilde \gamma$,
a normal geodesic from $\tilde p$ to $\gamma(\tilde p)$.
For each $\tilde \gamma(i)$, $1\le i\le |\gamma|$,  there is a pre-image of $p$,
say $\tilde p_i=\gamma_i(\tilde p)$, such that $|\tilde \gamma(i)\tilde p_i|\leq d$. Because $|\tilde p_i\tilde p_{i+1}|=|\gamma_i(\tilde p)\gamma_{i+1}(\tilde p)|\leq 2d+1$, we have that
$\gamma_i^{-1}\gamma_{i+1}\in \Gamma(2d+1)$. Then
$$\gamma=\gamma_1\gamma_1^{-1}\gamma_2\gamma_2^{-1}\cdots \gamma_{t-1}\gamma_{t-1}^{-1}\gamma_t=\gamma_1(\gamma_1^{-1}\gamma_2)\cdots (\gamma_{t-1}^{-1}\gamma_t),$$
where $\gamma_t=\gamma$, $t\leq R$ and $\gamma_1\in \Gamma(d+1)$. Thus $\gamma\in \Gamma(R, S)$.
\qed\enddemo

\vskip4mm

\subhead b. Groups with at most polynomial growth
\endsubhead

\vskip4mm

A group $N$ is nilpotent, if $N$ satisfying the following descending sequence,
$$N=N_0\vartriangleright N_1\vartriangleright\cdots \vartriangleright N_s=e,
\quad N_{i+1}=[N,N_i].$$
If $N$ is also finitely generated, then one can construct a symmetric generating
set, that adapts the descending sequence, as follows: starting with generators for $N_{s-1}$, $g_{s-1,1},...,g_{s-1,k_{s-1}}$, extend it to a generating set for $N_{s-2}$
by adding elements from $N_{s-2}-N_{s-1}$. Repeating this $(s-1)$-steps, we obtain
a generating set $S_0$.  Clearly, $S=S_0\cup S_0^{-1}$ is a finite symmetric generating set,
called a graded symmetric generating set.

Using a graded symmetric generating set, Wolf (\cite{Wo}) showed that $N$ has
at most polynomial growth; a property independent of $S$.

\proclaim{Theorem 1.3} Let $\Lambda$ be a finitely generated
nilpotent group. Then $\Lambda$ has at most
polynomial growth. Precisely, given any graded symmetric generating set, $\{\gamma_1, \dots, \gamma_m\}$,
there are polynomials $\{p_1,\dots,p_m\}$ such that for all $R\ge 1$ and any $\gamma\in\Lambda(R)$,
$\gamma$ has an expression,
$\gamma=\gamma_1^{a_1}\cdots\gamma_m^{a_m}$ such that $a_i\le p_i(R)$.
\endproclaim

We point it out that our proof of Theorem 0.8 follows the approach in the
proof of Theorem 1.3 in [Wo]; see discussion in Section 3.

A group $\Lambda$ is called virtually nilpotent, if $\Lambda$ contains a nilpotent subgroup $N$ of finite index. We now extend Theorem 1.3 to virtually nilpotent groups.

\proclaim{Lemma 1.4} The conclusion of Theorem 1.3 holds for any virtually nilpotent group $\Lambda$, i.e., there is a generator system $S=\{\tau_1,\dots, \tau_{m}, \beta_1,\dots, \beta_c\}$ such that for any $\gamma\in \Lambda(R,S)$, $\gamma=\tau_1^{a_1}\cdots\tau_m^{a_m}\beta_q$ for some $1\le q\le c$ such that $a_i\le p_i(R)$.
\endproclaim

We need some preparation.

\proclaim{Lemma 1.5} Assume that a group $\Gamma$ with a finite generating set $S$ satisfies an exact sequence, $1\to F\to \Gamma\overset{\pi}\to \longrightarrow \Lambda\to 1$.

\noindent \rm{(1.5.1)} If $F$ is
finite, then $h_w(\Gamma,S)=h_w(\Lambda,\bar S)$, where $\bar S=\pi(S)$. In particular,
$\Lambda$ has at most polynomial growth implies that $\Gamma$ does.

\noindent \rm{(1.5.2)} If $\Lambda$ is finite and $F\cap S$ generates $F$,
then for any $R>0$,
$$|\Gamma(R,S)|\leq |\Lambda|\cdot  |F(KR, F\cap S)|,$$
where $K>0$ is a constant (may depend on $S$).
\endproclaim

\demo{Proof} (1.5.1) follows from the fact that given any integer $R>0$,
$$|\Lambda(R,\bar S)|\le |\Gamma(R,S)|\le |F|\cdot
|\Lambda(R,\bar S)|.$$

(1.5.2) Since $\Lambda$ is finite, $F$ is also finitely generated. Let us take a symmetric generating set of $F$, $B_F=\{\alpha_1,..,\alpha_m\}$, and then add $B_0=\{\beta_1,...,\beta_c\}$ to form a set of generators for $\Gamma$,
$S=\{\alpha_1,..., \alpha_m,\beta_1,...,\beta_c\}$, where  $c=|\Lambda|$ and
$\Gamma$ is a disjoint union of costs,
$$\Gamma=\beta_1F,...,\beta_cF.$$
It suffices to prove (1.5.2) for $S$ and $B_F$.

Given any $\gamma=\gamma_1\gamma_2\cdots\gamma_R\in \Lambda_0(R,S)$ with each $\gamma_i\in S$,  we shall rewrite $\gamma$ so that all factors of $\gamma$ which are $\beta_i$'s move to far right via successive conjugations by $\beta_i$'s.
Then $\gamma$ can be written in the form $\gamma=\alpha\beta_{l},$
where $\beta_l\in B_0$ and $\alpha\in F$. If $|\alpha|_{B_F}\leq KR$ for some constant $K>0$,
then clearly (1.5.2) holds for $S$ and $S\cap F=B_F$.

Assume
$$\beta_i\beta_j=\alpha_{ij}\beta_{l(ij)}, \quad \beta_i\alpha_k\beta_i^{-1}=\alpha'_{ik},$$
where $1\leq i, j\leq c$, $1\leq k\leq m$, $\alpha_{ij}, \alpha'_{ij}\in F$, $\beta_{l(ij)}\in B_0$. Let
$$K=\max\{ |\alpha_{ij}|_{B_1},|\alpha'_{ik}|_{B_1},\, \beta_i\beta_j=\alpha_{ij}\beta_{l(ij)}, \beta_i\alpha_k\beta_i^{-1}=\alpha'_{ik}, 1\leq i, j\leq c, 1\leq k\leq m\}.$$

For $\gamma=\gamma_1\cdots\gamma_R$, if $\gamma_1=\beta_i$, then for $\gamma_2=\beta_j$,
$$\gamma=\alpha_{ij}\beta_{l(ij)}\gamma_3\cdots \gamma_R$$
and for $\gamma_2=\alpha_k$, $$\gamma=\alpha'_{ik}\beta_i\gamma_3\cdots\gamma_R.$$
Now consider $\gamma_3$ and repeat the process, and so do $\gamma_j$, $4\leq j\leq R$. And we have that $\gamma=\alpha\beta_l$, where $|\alpha|_{B_F}\leq KR$ and $\beta_l\in B_0$. \qed
\enddemo

\demo{Proof of Lemma 1.4}

Assume that $H$ is a nilpotent subgroup of $\Lambda$
such that $[\Lambda: H]=c>1$. Without loss of generality, we may
assume that $H$ is a normal subgroup of $\Lambda$. Assume $\Lambda=\bigcup_{i=1}^c \beta_iH$.

Let $S=B_H\cup B_0$ be a generator system of $\Lambda$ such that $\left<B_H\right>=H$ and $B_0=\{\beta_1,\beta_2,\dots \beta_c\}$. By the proof of (1.5.2), for each $\gamma\in \Lambda(R,S)$, $\gamma$ can be rewritten as
$$\gamma=\alpha \beta_l,$$
where $\alpha\in H(KR,B_H)$ and for some $1\leq l\leq c$.
Now by the fact $H$ is nilpotent and Theorem 1.3, we derive the result.
\qed\enddemo

\vskip4mm

\subhead c. Algebraic entropy of solvable groups
\endsubhead

\vskip4mm

A group $G$ is called solvable, if $G$ satisfies the following descending sequence,
$$G_0\vartriangleright G_1\vartriangleright\cdots \vartriangleright G_s=e,\quad G_0=G, G_{i+1}=[G_i,G_i].$$
By definition, a nilpotent group is solvable, but a solvable group may have a positive algebraic entropy.

\proclaim{Lemma 1.6 (\cite{Ti})}  Let $\Gamma$ be linear group over a field of characteristic $0$.

\noindent {\rm (1.6.1)} Either $\Gamma$ has a non-abelian free subgroup or $\Gamma$ possesses a solvable subgroup of finite index.

\noindent {\rm (1.6.2)} If $\Gamma$ is also finitely generated, then $h_w(\Gamma,S)=0$ implies that $\Gamma$ is virtually nilpotent (cf. \cite{Mi2},\cite{Wo}).
\endproclaim

\proclaim{Lemma 1.7 (\cite{Os1}, \cite{BH})} Let $\Gamma=\Bbb Z^k\rtimes_\rho\Bbb Z$ with $\rho: \Bbb Z
\to \op{Aut}(\Bbb Z^k)=\op{SL}(k,\Bbb Z)$. Then $h_w(\Gamma)=0$ if and only if $\Gamma$ is
virtually nilpotent, or equivalently, if and only if all eigenvalues of $\rho$ have norm one.
Furthermore, if $\lambda_{\max}(\rho)$ denotes the maximum for the norm
of eigenvalues of $\rho(1)$, then for any symmetric generating set $S$,
$$h_w(\Gamma,S)\ge \ln 2\cdot \frac {\ln \lambda_{\max}(\rho)}{\ln 2+5\ln \lambda_{\max}(\rho)}.$$
\endproclaim

\vskip4mm

\subhead d. Stabilities of fundamental groups in GH convergence
\endsubhead

\vskip4mm

The references for this subsections are [FY] (cf. [Ro]) and [KW].

Let $X_i$ be a sequence of compact length spaces, let
$\pi_i: (\tilde X_i,\tilde p_i)\to (X_i,p_i)$ be the
universal cover with the pullback length metrics, and let
$\Gamma_i=\pi_1(X_i,p_i)$.

Assume that $(X_i,p_i)@>GH>>(X,p)$ and $X$ is compact i.e., there is a sequence of maps, $f_i: (X_i,p_i)\to (X,p)$, called $\epsilon_i$-Gromov-Hausdorff approximations (briefly, GHA) with $\epsilon_i\to 0$, such that $f_i(p_i)=p$, $|d_X(f_i(x_i),f_i(x_i'))-d_{X_i}(x_i,x_i')|<\epsilon_i$, and
$f_i(X_i)$ is $\epsilon_i$-dense in $X$. Passing to a subsequence,
$(\tilde X_i,\tilde p_i,\Gamma_i)$ equivariant converges to
$(\tilde X,\tilde p,G)$ such that the following diagram commutes:
$$\CD (\tilde X_i,\tilde p_i,\Gamma_i)
@> (\tilde f_i,\phi_i)>> (\tilde X,\tilde p,G)\\
@ VV \pi_i V@VV \pi  V\\
(X_i,p_i) @> f_i >> (X,p),\endCD\tag 1.8$$
where $G$ is a closed subgroup of isometries on $\tilde X$, i.e.,
there is a sequence of pair of maps, $(\tilde f_i,\phi_i): (\tilde X_i,\tilde p_i,\Gamma_i)\to (\tilde X,\tilde p,G)$, which are $\epsilon_i$-equivariant GHA $(\tilde f_i,\phi_i)$ consists of two $\epsilon_i$-GHAs $\tilde f_i: (\tilde X_i,\tilde p_i)\to
(\tilde X,\tilde p)$, $\phi_i: \Gamma_i\to G$ ($\Gamma_i$
and $G$ are equipped with the metrics induced from $\tilde X_i$ and $\tilde X$) such that the two actions almost commutes in terms of $\tilde f_i$ and $\phi_i$:
$$d_X(\tilde f_i(\gamma_i(\tilde x_i)),\phi_i(\gamma_i)(\tilde f_i(\tilde x_i))
<\epsilon_i,\quad \tilde x_i\in B_{\epsilon_i^{-1}}(\tilde p_i)\subset \tilde X_i, \, \gamma_i\in \Gamma_i(\epsilon_i^{-1}).$$
(for details, see \cite{FY}, cf. \cite{Ro}).

A stability question on (1.8) is: to what extend, does $\Gamma_i$
become ``stable''?

\proclaim{Theorem 1.9 (\cite {FY})} Let $X_i$ be compact length spaces
in (1.8) such that $X$ is compact and $G/G_0$ is discrete, where $G_0$ is the identity component of $G$. Then there exist $\epsilon=\epsilon(X)>0$ and $R_i\to \infty$ such that the
subgroup $\Gamma_{i,\epsilon}$ generated by $\{\gamma_i\in \Gamma_i,\,\, |\tilde q_i\gamma_i(\tilde q_i)|<\epsilon, \,\,\tilde q_i\in B_{R_i}(\tilde p_i)\}$
is normal in $\Gamma_i$ and the equivariant GH-convergence induces an isomorphism: $\phi_i: \Gamma_i/\Gamma_{i,\epsilon}\to G/G_0$.
\endproclaim

Note that if $X_i=M_i$, a compact Riemannian manifold of $\op{Ric}_{M_i}\ge -(n-1)$, then $\Gamma_{i,\epsilon}$ is generated by $\{\gamma_i\in \Gamma_i,\,\, |\tilde q_i\gamma_i(\tilde q_i)|<\epsilon, \,\,\tilde q_i\in B_1(\tilde p_i)\}$.

We now specify $X_i=M_i$, a compact Riemannian $n$-manifold such that
$$\op{Ric}_{M_i}\ge -(n-1),\quad d\ge \op{diam}(M_i).$$
By \cite{CN}, $G$ is a Lie group, and thus $G/G_0$ is discrete.

\proclaim{Lemma 1.10}\rm{(\cite{PR})} Let $M_i@>GH>>X$ be a sequence of
compact $n$-manifolds as in the above which satisfies (1.8). If there is
$v>0$ such that $\op{vol}(M_i)\ge v$,
then for $i$ large, the isomorphism class of $\pi_1(M_i)$ is independent of $i$.
\endproclaim

Note that Lemma 1.10 improves a result of \cite{An} on the finiteness of the fundamental
groups of $n$-manifolds. The proof of Lemma 1.10 is to show that $\Gamma_{i,\epsilon}=e$.

According to \cite{KW}, $\Gamma_{i,\epsilon}$ contains a nilpotent subgroup
$N_i$ which is normal in $\Gamma_i$ and whose index is bounded by a constant depending on $n$ and $d$. The following stability property on $(\Gamma_i,N_i)$
in \cite{KW} will serve as a starting point for our proof of Theorem 0.3
(see Theorem 1.12 below):

\noindent (1.11.1) The torsion group of
$N_i$, $\op{Tor}(N_i)$, is normal in $\Gamma_i$ and, after modulo $\operatorname{Tor}(N_i)$
the length of the nilpotency of $N_i$, $k_i\le k$ ($\le n-1$):
$$N_i=N_{i1}\vartriangleright\cdots \vartriangleright N_{ik_i}\vartriangleright N_{i(k_i+1)}=\op{Tor}(N_i), \qquad [N,N_h]\subset N_{h+1}.$$

\noindent (1.11.2) Passing to a subsequence, we may assume for $1\le h\le k$, $N_{i,h}/N_{i,h+1}
\cong \Bbb Z^{m_h}$,  and  the homomorphism, $\bar \rho_h: \Gamma_i/N_i
\to \text{Aut}(N_{i,h}/N_{i,h+1})$, are independent of $i$ (see discussion before Theorem 0.7).

\proclaim{Theorem 1.12 (\cite{KW})}  (1.12.1) Given $n, d>0$,  there are finitely many groups, $F_1,
\dots, F_k$, each of which is finitely presented in terms of a fixed generator system $S_1,\dots,S_k$ respectively  (according to Lemma 9.2 (b) in [KW]), such that the following holds:

For any compact $n$-manifold $M$ with $\op{Ric}_M\geq -(n-1)$ and $\op{diam}(M)\leq d$ there is a nilpotent normal
subgroup $N\vartriangleleft \pi_1(M)$ which has a nilpotent basis of length $\leq n-1$ such that $\pi_1(M)/N$ is isomorphic to $F_j$, where $S_k$ corresponds to a finitely many loops in $\pi_1(M)$ at a base point whose length and number is bounded by a constant $c(n,d)$.

\noindent (1.12.2) In addition to (1.12.1), one can choose a finite collection of irreducible rational representations
$\rho_i^j:F_i\to GL(n_i^j,\Bbb Q)$ $(j=1,\cdots, \mu_i, i=1,\cdots, k)$ such that for a suitable choice of the isomorphism $\pi_1(M)/N \cong F_i$ the following holds:

There is a chain of subgroups $N=N_{1}\vartriangleright \cdots \vartriangleright N_{h_0}\vartriangleright N_{h_0+1}=\op{Tor}(N)$
which are all normal in $\pi_1(M)$
such that for any $1\le h\le h_0$, $[N,N_h]\subset N_{h+1}$ and $N_h/{N_{h+1}}$ is free abelian. Moreover, the action of $\pi_1(M)$ on $N$ by
conjugation induces an action of $F_i$ on $N_h/{N_{h+1}}$ and the induced representation $\bar \rho_hF_i\to \op{GL}(N_h/{N_{h+1}}
\otimes_{\Bbb Z}\Bbb Q)$ is isomorphic to $\rho_i^j$ for a suitable $j=j(h), h=0,\cdots, h_0-1$.
\endproclaim


	
\vskip4mm

\head 2. Proofs of Theorem 0.3 and 0.8 by Assuming Theorem 0.7
\endhead

\vskip4mm

As explained in the introduction, our approach is to establish, using that $h(M)$ is small, (0.7.1) and (0.7.2) (see (2.1.1) and (2.1.2) below), and
apply Theorem 0.7 to conclude Theorem 0.3 and Theorem 0.8.

\proclaim{Lemma 2.1} Given $n ,d>0$, there exists $\epsilon(n, d)>0$ such that if a compact $n$-manifold $M$ satisfies that
$$\op{Ric}_{M}\ge -(n-1),\quad \op{diam}(M)\le  d, \quad  h(M)\leq \epsilon(n, d),$$
and $(\Gamma,N)$ satisfies Theorem 1.12, then the following properties hold:

\noindent (2.1.1) $h_w(\Lambda, \pi(\Gamma(2d+1)))=0$, where $\pi:\Gamma\to \Lambda$ is the quotient map.

\noindent (2.1.2) For any $\gamma\in \Gamma$, the eigenvalues of $\bar \rho_h(\bar \gamma)$ have norm one, provided $\bar \gamma=\pi(\gamma)$ has an infinite order (equivalently, $\left<\gamma, N_{h}\right>/N_{h+1}$ is virtually nilpotent).
\endproclaim

\demo{Proof} By the Gromov's compactness, it is equivalent to consider
a sequence $M_i$ in (0.5), and show that (2.1.1) or (2.1.2) hold for large $i$.
From the discussion in Section 1, we may assume an equivariant convergent sequence $(\tilde M_i,\tilde p_i,\Gamma_i)
@> \op{eqGH} >> (\tilde X,\tilde p,G)$ as in (1.8), and there is $\epsilon>0$ such that
$i$ large, there is an isomorphism, $\bar \psi_i: \bar \Gamma_i = \Gamma_i/\Gamma_{i,\epsilon}\to \bar G=G/G_0$, which is also an $\Psi(\epsilon_i)$-GHA with respect to the natural metrics on $\Gamma_i$ and $G$ induced from $\tilde M_i$ and $\tilde X$.

(2.1.1) Note that there is a natural surjective homomorphism, $\Lambda_i \to \bar \Gamma_i$, $\gamma N_i\mapsto \gamma \Gamma_{i,\epsilon}$, with a finite kernel. By (1.5.1), it suffices to show
$h_w(\bar \Gamma_i,\bar \Gamma_i(2d+1))=0$.
Note that the discreteness of $\bar G(2d+1)$ implies that $\bar G(2d+1)=\bar \psi_i(\bar \Gamma_i(2d+1))$ is independent of large $i$. It follows that
$$h_w(\bar G,\bar G(2d+1))=h_w(\bar \Gamma_i, \bar \Gamma_i(2d+1))\le h_w(\Gamma_i,\Gamma_i(2d+1))\to 0,$$
and thus $h_w(\bar G,\bar G(2d+1))=0$.

(2.1.2) If (2.1.2) fails, then we may assume a sequence, $\bar \gamma_i\in \Lambda_i$, of infinite order such that $\bar \rho_h(\bar \gamma_i)$ has an eigenvalue of norm $\ne 1$. We may assume the maximum of the norm of
eigenvalues of $\bar \rho_h(\bar \gamma_i)$ is $>1$.

By Theorem 1.12, there is a large $i_0$ such that for all $i\ge i_0$,
$\bar \phi_i: \Lambda_i\to \Lambda\subset G/G_0$ is an isomorphism, and the induced
$\bar \rho_{i,h}: \Lambda_i\to \op{Aut}(\Bbb Z^{m_h})$ is independent
of $i$. For all $i\ge i_0$, let $\bar \beta_i=\bar \phi_i^{-1}(\bar \phi_{i_0}(\gamma_{i_0}))$. Then one can choose $\beta_i\in \bar \psi_i^{-1}(\bar \beta_i)$ such that
the maximal displacement of $\beta_i$ on $B_1(\tilde p_i)$ is bounded above by
two times of that of $\gamma_{i_0}$, say $L$.

Let $\bar \eta_i: \Gamma_i\to \Gamma_i/N_{i,h+1}$ be the projection, and let $\Bbb Z^{m_h}=N_{i,h}/N_{i,h+1}$. Because the length of $\beta_i$ is bounded by $L$, passing to a subsequence we may assume that the subgroup of $\Gamma_i/N_{i,h+1}$,
$H=\left<\bar \beta_i, \Bbb Z^{m_h}\right>\cong\Bbb Z^{m_h}\rtimes_{\bar \rho_h(\bar \beta_i)}
\Bbb Z$ is independent of large $i$. Let $\bar S$ be the set of generators formed by a canonical basis of $\Bbb Z^{m_h}$ and $\bar{\beta_i}$, and let $S_i\subset \eta_i^{-1}(\bar S)$ be of 1-1 correspondence to $\bar S_i$ whose elements are of length $\le L$.
By Lemma 1.7, $h_w(H,\bar S)\ge
\ln 2\cdot \frac {\ln \lambda_{\max}(\bar \beta_i)}{\ln 2+5\ln \lambda_{\max}(\bar \beta_i)}>0$, a contradiction to the fact that
$$h_w(H,\bar S)\le h_w(\eta_i^{-1}(H), S_i)\le h_w(\Gamma_i,\Gamma_i(L))\le L\cdot h(M_i)\to 0.$$
\qed\enddemo

\demo{Proof of Theorem 0.3 by assuming Theorem 0.7}

Let $M$ satisfy the conditions of Theorem 0.3. By \cite{Gr}, it reduces to show that
$\Gamma$ has at most a polynomial growth. By Lemma 2.1, $(\Gamma,N)$
satisfies (0.7.2), and $h_w(\Lambda,S)=0$ which implies that (0.7.1) with a
positive answer to the Milnor Problem 0.1. By now Theorem 0.3 follows from
Theorem 0.7
\qed\enddemo

\demo{Proof of Theorem 0.8 by assuming Theorem 0.7}

Let $M$ be as in Theorem 0.8. By \cite{Gr}, it reduces to show that
$\Gamma$ has at most polynomial growth. By Lemma 2.1, $(\Gamma,N)$
satisfies (0.7.2) and $h_w(\Lambda,S)=0$.

(0.8.1) First, $\Gamma$ and $N$ (nilpotent) are elementary amendable.
Then $\Lambda=\Gamma/N$ is also elementary amenable. Because $h_w(\Lambda,S)=0$,
$\Lambda$ is virtually nilpotent (\cite{Ch}). We then apply Theorem 0.7
to conclude that $\Gamma$ has a polynomial growth.

(0.8.2) We claim that $\Gamma$ has a solvable subgroup of a finite index,
and in particular $\Gamma$ is elementary amendable, and (0.8.2) reduces to
(0.8.1).

By Lemma 1.6, the claim follows if we show that any two elements $\alpha,\beta$ of $\pi_1(M)$, the subgroup $\left<\alpha,\beta\right>$ is not a non-abelian free group.

Arguing by contradiction, let us assume that $\left<\alpha,\beta\right>$ is non-abelian free while $h(M)<\epsilon(n,d)$, which is given in Lemma 2.1. Let $\bar\alpha, \bar\beta$ be the projection of $\alpha,\beta$ in $\Lambda$. Since $h_w(\left<\bar\alpha,\bar \beta\right>)=0$, there is a non-trivial relation in $\Lambda$, say $$\bar\alpha^{\pm i_1}\bar\beta^{\pm i_2}\cdots \bar\alpha^{\pm i_1}\bar\beta^{\pm i_m}=\bar e.$$
That is, $\gamma=\alpha^{\pm i_1}\beta^{\pm i_2}\cdots \alpha^{\pm i_m}\beta^{\pm i_m}\in N$. Because $\left<\alpha,\beta\right>$ is free, $\gamma\neq e$. Without loss of generality, we assume $\gamma\neq \beta$. Because any subgroup of a free group is also free (the Nielsen-Schreier theorem), $\left<\gamma,\beta\right>$ is a free group by two generators. On the other hand, $\left<\gamma,\beta\right>$ is a subgroup of $\left<\beta, N\right>$ which by Lemma 2.1 is virtually nilpotent. We derive
$$2=h_w(\left<\gamma,\beta\right>, \{\gamma,\beta\})\le h_w(\left<\beta, N\right>)= 0,$$ a contradiction.
\qed\enddemo

\vskip4mm

\head 3. Outline of Proof of Theorem 0.7.
\endhead

\vskip4mm

Our approach to Theorem 0.7 is motivated by the proof of Theorem 1.3 (Theorem 3.2 in \cite{Wo})
via the method of `successively conjugate operations' (see (3.3) and (3.6)). The approach is also readily to apply for Theorem 0.9 in a more general setting where $\Lambda$ is not virtually nilpotent. Because
the proof is tedious and technical, for convenience of readers, we provide an outline.

Let $1\to N\to \Gamma \to \Lambda\to 1$ be the exact sequence in  Theorem 0.7. Without loss of generality, we may assume that $N$ is torsion free and $\Lambda$ is nilpotent. Indeed,
since $N$ is normal, $\op{Tor}(N)$ is normal in $\Gamma$.
If $N'=N/\op{Tor(N)}$ is non-trivial, then let us consider the exact sequence
$$1\to N'\to \Gamma' \overset{\pi'}\to{\to} \Lambda\to 1.$$
Since the conjugate actions $\bar\rho_h$ of $\Lambda$ on $N_h'/N_{h+1}'=N_h/N_{h+1}$ are the same as the original action of $\Lambda$ on $N_h/N_{h+1}$, (0.7.2) holds for the new exact sequence.

\proclaim{Lemma 3.1} Let $N$ be a finitely generated torsion nilpotent group.
Then $N$ is a finite group.
\endproclaim

\demo{Proof} Consider the descending sequence of normal subgroups,
$$N\vartriangleright N_1\vartriangleright\cdots \vartriangleright N_{k+1}=e,\quad N_{i+1}=[N,N_i].$$
We will proceed by induction on $k$. For $k=0$, $N$ is finite because $N$ abelian. Assume Lemma 3.1 holds for $k$. Consider $N\vartriangleright N_1$, $|N|\le |N/N_1|\cdot |N_1|$. By induction, $N_1$ is finite. Because $N/N_1$ is abelian and thus finite.
\qed\enddemo

By Lemma 3.1, $\op{Tor}(N)$ is finite, and thus by Lemma 1.5 algebraic entropy and polynomial growth of subgroups in $\Gamma'$ coincides with the corresponding ones in $\Gamma$. Since $\Lambda$ has a nilpotent subgroup $\Lambda_0$ of finite index, $\Gamma_0'=(\pi')^{-1}(\Lambda_0)$ is of finite index in $\Gamma'$. By Lemma 1.5 again, the growth of $\Gamma_0'$ is the same as $\Gamma$.

As the first step, we will prove that the generator system $S=B_1\cup S_0$ of $\Gamma$, where $B_1$ is a finite graded basis for $N$ and $S_0=\left\{\gamma_1,\dots,\gamma_m\right\}$ be given in (0.7.2) satisfies the following property.

\proclaim{Lemma 3.2}
Any word representation,
$w_\gamma=v_1\cdots v_R$, in terms of $S$ can be rewritten into a product of a form,
$$w_\gamma=(\alpha_{i_1}\cdots \alpha_{i_s})\cdot (\gamma_{j_1}\cdots \gamma_{j_t}), \tag{3.3}$$
where either $\alpha_{i_l}$ or $\alpha_{i_l}^{-1}\in B_1$ and $\gamma_j\in S_0$, $t\le R$, such that
$$ s\le P_0(R), \qquad\text{where $P_0$ is a polynomial.} \tag 3.4$$
In particular, the ``horizontal'' lifting subgroup $\Gamma_0=\left<S_0\right>$ of $\Lambda$ satisfies $h_w(\Gamma_0,S_0)=h_w(\Gamma,S)$.
\endproclaim

Secondly, by a similar idea of the proof of Theorem 1.3, we will prove that $\gamma_{j_1}\cdots \gamma_{j_t}$ in (3.3) can be written into a standard form in terms of a graded symmetric basis $\bar T=\{\bar\tau_1,\cdots,\bar\tau_{m_l}\}$ of $\Lambda$ (see Theorem 1.3):

\proclaim{Lemma 3.5}
After fixing a subset $T=\{\tau_1,\dots,\tau_{m_l}\}$ of $\Gamma$ whose projection $\bar T$ is a graded symmetric basis of $\Lambda$, any $\gamma\in \Gamma_0=\left<S_0\right>$ of length $R$ can be written as
$$\gamma=
\alpha\tau_{m_l}^{P_{m_l}(\gamma)}\cdots\tau_1^{P_1(\gamma)}, \quad |\alpha|_{B_1}\leq P(R), \tag 3.6$$
where  $P_i(\gamma)$ are bounded by polynomials $P_i(R)$ in $R$.
\endproclaim

Because $N$ is of polynomial growth, it is clear that Theorem 0.7 follows from Lemma 3.2 and Lemma 3.5.

The main task here is to show (3.3) and (3.4), whose idea is given as follows. Let us choose a finite graded basis for $N$: for $1\le h\le k$, extending the gradated basis to a symmetric set generating
set for $N_h$
$$B_h=\{\alpha_{n_{h-1}+1},\dots, \alpha_{n_h}; \alpha_{n_h+1},\dots,\alpha_{n_{h+1}};\dots;\alpha_{n_{k-1}+1},
\dots, \alpha_{n_k}\},\tag{3.7}$$
such that $\bar B_h=\{\bar{\alpha}_{n_{h-1}+1},
\cdots, \bar{\alpha}_{n_h}\}$ is a basis of $N_h/{N_{h+1}}$, where
$\bar B_h$ denote the projection of $B_h$ into $N_h/N_{h+1}$. Let $S_0=\{\gamma_1,...,\gamma_m\}\subset \Gamma$ be given as in (0.7.2) such that its projection $\bar S_0=\{\bar \gamma_1,...,\bar \gamma_m\}$ is a graded symmetric
generating set for $\Lambda=\Gamma/N$. Then for all $1\le h\le k_1$, $\{\bar \rho_h(\bar \gamma_1),...,
\bar\rho_h(\bar \gamma_m)\}$ contains a graded symmetric generating set for $\op{Im}(\bar \rho)$,
where $\bar \rho_h: \Lambda\to \op{Aut}(N_h/N_{h+1})$ is induced by the conjugate action, $\bar\rho_h(\bar \gamma)(\alpha N_{h+1})=\gamma \alpha \gamma^{-1} N_{h+1}$. Then $S=B_1\cup S_0$ forms a generator system of $\Gamma$.

For any $\gamma\in \Gamma_w(R,S)$ with a word representation,
$w_\gamma=v_1\cdots v_R$, such that either $v_i$ or
$v_i^{-1}\in S$, we can achieve (3.3) via the following successive conjugate operations on $w_\gamma$:

Let $\gamma_{j_1}$ denote the first element from left that is in $S_0$, and let
$v_q$ denote the right adjacent element. If $v_q$ or $v_q^{-1}=\alpha_p\in B_1$, then
$$\gamma_{j_1}\alpha_p=\gamma_{j_1}\alpha_p\gamma_{j_1}^{-1}\cdot \gamma_{j_1}=\rho(\gamma_{j_1})
(\alpha_p)\cdot \gamma_{j_1}, \quad \rho(\gamma_{j_1})(\alpha_p)\in N_h, \quad n_{h-1}< p
\le n_h.$$
If $v_q=\gamma_{j_2}\in S_0$, then we consider the next right adjacent element and iterate
the above process. In such a way, we will perform the conjugate operations at most $R$ times, each one
looks like,
$$\split \gamma_{j_1}\cdots \gamma_{j_t}\alpha_p&=(\gamma_{j_1}\cdots
\gamma_{j_t}\alpha_p\gamma_{j_t}^{-1}\cdots \gamma_{j_1}^{-1})\cdot
(\gamma_{j_1}\cdots \gamma_{j_t})\\&=\rho(\gamma_{j_1}\cdots \gamma_{j_t})
(\alpha_p)\cdot (\gamma_{j_1}\cdots \gamma_{j_t}),\endsplit$$
and $\rho(\gamma_{j_1}\cdots \gamma_{j_t})(\alpha_p)\in N_{h-1}$, $n_{h-1}<p\le n_h$. We claim
that for each $1\le h\le k$,
$$|\rho(\gamma_{j_1}\cdots \gamma_{j_t})(\alpha_p)|_{B_1}\le P(R),\tag 3.8$$
where $P$ is a fixed polynomial in $R$ (independent of the word representation $w_\gamma$). From the above
discussion, (3.8) implies (3.4), $s\le (R-1)\cdot P(R)$.

We will achieve (3.8) by reverse induction on $h$, starting with $h=k$
i.e., $N_k$ is a free abelian group. Since $
\op{Im}(\bar \rho_k)$ is nilpotent, we can have the
equation (Theorem 1.3),
$$\rho(\gamma_{j_1}\cdots \gamma_{j_t})(\alpha)=\bar \rho_k(\gamma_1^{l_1}
\cdots \gamma_m^{l_m})(\alpha),\quad \alpha\in N_k,$$
and $l_j$ is at most polynomial in $t$.

The assumption that for any $\gamma_j\in S_0$, all eigenvalues of $\bar \rho_k(\bar \gamma_j)$
have norm one implies that $\left<\gamma_j,N_k\right>$ is virtually nilpotent. Because all
eigenvalues (possibly complex) of $\rho(\gamma_j)$ have norm one, if $A(\gamma_j)$ denotes
a matrix representation of $\rho(\gamma_j)$, then there is a complex matrix $X$
such that $X^{-1}A(\gamma_j)X$ is a Jordan matrix with diagonal entires norm one.
Consequently, the entries of $(A(\gamma_j))^{l_j}$ is bounded above by a polynomial
which also depends on $\gamma_j$.
Hence, we may assume that $$\gamma_j^{l_j}\alpha_p = \alpha_{n_{k-1}+1}
^{P^j_{p,n_{k-1}+1}(l_j)}\cdots \alpha_{n_k}^{P^j_{p,n_k}(l_j)}\gamma_j^{l_j},$$
where $P_{i, u}^j(l_j)$ are polynomials in $l_j$. It follows that for $\alpha_p\in B_k$,
$$\gamma_{j_1}\cdots \gamma_{j_t}\alpha_p\gamma_{j_t}^{-1}\cdots \gamma_{j_1}^{-1}=
\alpha_{n_{k-1}+1}^{P^1_{p,n_{k-1}+1}(l_1)+\cdots +P^m_{p,n_{k-1}+1}(l_m)}
\cdots \alpha_{n_k}^{P^1_{p,n_k}(l_1)+\cdots +P^m_{p,n_k}(l_m)},$$
which finish the proof of (3.8) for $h=k$.

The general case of (3.8) where $N_h$ ($h<k$) is not abelian will be treated in the next section.
\vskip4mm

\head 4. Proof of Theorem 0.7
\endhead

\vskip4mm

The proof of Theorem 0.7 will follow the approach and setup in the last section. It suffices to prove Lemma 3.2 and Lemma 3.5.

Recall that we have prove (3.8) for $N_k$ in Section 3, as formulated in the following.
\proclaim{Lemma 4.1} Let the assumptions be as in Theorem 0.7, where $N=\left<\alpha_1,\dots,\alpha_n\right>$ is abelian. Let $\Gamma_0=\left<S_0\right>$. Then for any $R>0$, $\gamma\in \Gamma_0(R,S_0)$, any $q\in \{1, \cdots, n\}$
$$\gamma\alpha_q \gamma^{-1}=\alpha_1^{P_{q1}(\gamma)}\cdots \alpha_{n}^{P_{qn}(\gamma)},$$
where $P_{qi}(\gamma)$ are bounded by polynomials $P_{qi}(R)$ in $R$, for $i=1,\cdots, n$.
\endproclaim

Now let us make some preparation for the proof of (3.8) for $N_h$ ($h<k$). Since $N_h$ is not abelian, we may not have a matrix presentation for $\rho(\gamma_j)$,
which may cause some complexity during counting the length of words. Because of this,
we shall first fix word expression
for each $\rho(\gamma_j)(\alpha_i)$, $\gamma_j\in S_0$ and $\alpha_i\in B_h\setminus B_{h+1}$ $(h=1,\dots, k, B_{k+1}=e)$:
$$\rho(\gamma_j)(\alpha_i)=\alpha_{n_{h-1}+1}^{a_{i,n_{h-1}+1}(\gamma_j)}\cdots
\alpha_{n_h}^{a_{i,n_h}(\gamma_j)}\cdots \alpha_{n_k}^{a_{i,n_k}(\gamma_j)}.$$
For $\bar \rho_h: \Lambda\to \op{Aut}(N_h/N_{h+1})$, similarly we can express
$$\bar \rho_h(\bar \gamma_j)(\bar \alpha_i)=\bar \alpha_{n_{h-1}+1}^{\bar a_{i,n_{h-1}+1}
	(\bar \gamma_j)}\cdots \bar \alpha_{n_h}^{\bar a_{i,n_h}(\bar \gamma_j)}.$$
Because the following diagram commutes for all $\gamma\in \Gamma$,
$$\CD N_h @> \rho(\gamma)>> N_h\\
@VV V @VV V\\N_h/N_{h+1}@>\bar \rho_h(\bar \gamma)>>N_h/N_{h+1},\endCD$$
and because the projection of $B_h$ on $\bar B_h$ is a basis for $N_h/N_{h+1}$, we conclude that
$$a_{i,p}(\gamma_j)=\bar a_{i,p}(\bar \gamma_j),\quad n_{h-1}+1\le p\le n_h. \tag 4.2$$
Let $$\beta_{i,h}(\gamma_j)=\alpha_{n_h+1}^{a_{i,n_h+1}(\gamma_j)}\cdots \alpha_{n_k}^{a_{i,n_k}(\gamma_j)}.$$
For each $\gamma=\gamma_{j_2}\gamma_{j_1}\in \Gamma$,
$$\aligned &\rho(\gamma_{j_{2}}\gamma_{j_{1}})(\alpha_i)
=\rho(\gamma_{j_2})\left(\alpha_{n_{h-1}+1}^{a_{i,n_{h-1}+1}(\gamma_{j_1})}\cdots
\alpha_{n_h}^{a_{i,n_h}(\gamma_{j_1})}\cdot \beta_{i,h}(\gamma_{j_1})\right)\\
=&\left(\rho(\gamma_{j_2})(\alpha_{n_{h-1}+1})\right)^{a_{i,n_{h-1}+1}(\gamma_{j_1})}
\cdots \left(\rho(\gamma_{j_2})(\alpha_{n_h})\right)^{a_{i,n_h}(\gamma_{j_1})}
\cdot
\rho(\gamma_{j_2})(\beta_{i,h}(\gamma_{j_1}))
\endaligned$$
By Theorem 1.3, the subword before $\rho(\gamma_{j_2})(\beta_{i,h}(\gamma_{j_1}))$ can be rewritten as
$$\aligned
&\left(\rho(\gamma_{j_2})(\alpha_{n_{h-1}+1})\right)^{a_{i,n_{h-1}+1}(\gamma_{j_1})}
\cdots \left(\rho(\gamma_{j_2})(\alpha_{n_h})\right)^{a_{i,n_h}(\gamma_{j_1})}\\
=&\alpha_{n_{h-1}+1}^{a_{i,n_{h-1}+1}(\gamma_{j_2}\gamma_{j_1})}\cdots
\alpha_{n_h}^{a_{i,n_h}(\gamma_{j_2}\gamma_{j_1})}\cdot \mu_{i,h}(\gamma_{j_2}
\gamma_{j_1}),\endaligned$$
where $\mu_{i,h}(\gamma_{j_2}\gamma_{j_1})$ denotes the remaining term that consists of $\alpha_{n_h+1},\dots,\alpha_{n_k}$.
Let $A=\max \{|a_{i,p}(\gamma_j)|, \,\, 1\le i,p\le n_k; 1\le j\le m\}$.
Then by Theorem 1.3, $\mu_{i,h}(\gamma_{j_2}\gamma_{j_1})$ has word length $\le P_0(n_k^2A^2)$.

Repeating the above we see that ($t\le R$)
$$\aligned\rho(\gamma_{j_t}\cdots\gamma_{j_{1}})(\alpha_i)
=&\alpha_{n_{h-1}+1}^{a_{i,n_{h-1}+1}(\gamma_{j_t}\cdots\gamma_{j_1})}\cdots
\alpha_{n_h}^{a_{i,n_h}(\gamma_{j_t}\cdots\gamma_{j_1})}\cdot\\
&\mu_{i,h}(\gamma_{j_t}\cdots\gamma_{j_1}) \cdot \rho(\gamma_{j_t})\left(\mu_{i,h}
(\gamma_{j_{t-1}}\cdots\gamma_{j_2}\gamma_{j_1})\right)\cdots\\
&\rho(\gamma_{j_t}\cdots \gamma_{j_{u-1}}) \left(\mu_{i,h}(\gamma_{j_u}\cdots
\gamma_{j_1})\right) \cdots \rho(\gamma_{j_t}\cdots\gamma_{j_2})\left(\beta_{i,h}
(\gamma_{j_{1}})\right).\endaligned$$

By the same reasoning for the identity in (4.2), from the above expression
we derive similar identities:
$$a_{i,p}(\gamma_{j_{t}}\cdots\gamma_{j_1})=\bar a_{i,p}(\bar \gamma_{j_{t}}\cdots\bar \gamma_{j_1}).\tag 4.3$$
By Lemma 4.1, we may assume that
$$\max\{|\bar a_{i,p}(\bar \gamma_{j_{t}}\cdots\bar \gamma_{j_1})|,\,\, n_{h-1}+1\le i,p\le n_h, \, 1\le j\le m\}\le \bar P_h(R).$$
With the above preparation, we are ready for

\demo{Proof of Lemma 3.2}

By the discussion in the outline of the proof in Section 3, the proof reduces to establish (3.8) for any $h$,
$1\le h\le k$. We shall proceed by inverse induction on $h$, where the case of $h=k$ follows from Lemma 4.1

Assume (3.8) for $h+1$ i.e., for any $\gamma=\gamma_{j_t}\cdots\gamma_{j_1}\in \Gamma$
of word length $t\le R$,
$$|\rho(\gamma_{j_t}\cdots\gamma_{j_{1}})(\alpha_i)|_{B_1}\le P_{h+1}(R),\quad \text{for any } \alpha_i\in
B_{h+1}.$$
It remains to check that (3.8) for $h$ i.e., there is a polynomial $P_h(R)$ such that
$$|\rho(\gamma_{j_t}\cdots\gamma_{j_{1}})(\alpha_i)|_{B_1}\le P_h(R),\quad \text{for any } \alpha_i\in
B_h.$$
For $n_{h-1}+1\le i,p\le n_h$, by (4.3) we are able to
apply Theorem 1.3 to conclude that $\mu_{i,h}(\gamma_{j_u}\cdots\gamma_{j_1})$ has word length $\le P_0(\hat R)$, where $\hat R=n_k^2\cdot A\cdot \bar P_h(R)$. Together with the inductive assumption, we derive
$$|\rho(\gamma_{j_t}\cdots \gamma_{j_{u-1}}) \left(\mu_{i,h}(\gamma_{j_u}\cdots\gamma_{j_1})
\right)|_{B_1}\le P_{h+1}(R)P_0(\hat R).$$
Then $$|\rho(\gamma_{j_t}\cdots\gamma_{j_{1}})(\alpha_i)|_{B_1}\le n_h\cdot \bar P_{h}(R)
+R\cdot P_{h+1}(R)\cdot P_0(\hat R)=P_h(R).$$
By now we have proved (3.8), which implies (3.4).\qed\enddemo

\remark{Remark 4.4}
It can be seen from its proof that,
Lemma 3.2 still holds after (0.7.1-2) are weakened to the condition that, instead of $\Lambda$ itself, one assumes that for each $1\le h\le k$, $\op{Im}(\bar\rho_h)$ is virtually nilpotent, and the image $\{\bar \rho_h(\pi(\gamma_i))\}$ contains a graded symmetric generating set for $\op{Im}(\rho_h)$ such that if $\pi(\gamma_i)$ has an infinite order, then eigenvalues of $\bar \rho_h(\pi(\gamma_i))$ have norm one. Indeed, by Lemma 1.4
$$\bar \rho_h(\gamma_{j_1}\cdots\gamma_{j_t})=\bar\rho_h(\gamma_1^{l_1}\cdots\gamma_m^{l_m}),$$
where $l_j$ is at most polynomial in $t$. Then (3.8) follows from the same argument.
\endremark

We now prove Lemma 3.5. As discussed at the beginning of Section 3, we have assumed that $\Lambda$ is nilpotent, and the generator system $S=B_1\cup S_0$, where $\bar S_0\subset \Lambda$ is a graded symmetric generator system of $\Lambda$ associated to
$$\Lambda_0=\Lambda\vartriangleright \Lambda_1\vartriangleright\cdots \vartriangleright \Lambda_{l+1}=e,\quad \Lambda_{i+1}=[\Lambda,\Lambda_i].$$
Let $\bar S_{p}=\{\bar \tau_{m_{p-1}+1}, \cdots, \bar \tau_{m_l}\}$ be the graded symmetric generator system of $\Lambda_p$.
Then any combination of elements in $S_0$ satisfies (3.8).

\demo{Proof of Lemma 3.5}

Let us consider the exact sequence $1\to N\to \Gamma_l\to \Lambda_l\to 1$, where $\Gamma_l=\pi^{-1}(\Lambda_l)$. Then by Lemma 3.2, any word $\gamma=v_1\cdots v_R$ in $\Gamma_l$ can be expressed as
$$v_1\cdots v_R=\alpha \cdot \tau_{j_1}\cdots \tau_{j_t}, \quad \text{where $|\alpha|_{B_1}\le P(R)$ and $\tau_{j_i}\in T_l.$}$$
Because $\Lambda_l$ is abelian, let us assume that for $m_{l-1}+1\le i,j\le m_l$,
$$\tau_i\tau_j=\alpha_{ij}\tau_j\tau_i, $$
and let $K=\max_{i,j,k}\{|\alpha_{ij}|_{B_1}, |\tau_i\alpha_k\tau_i^{-1}|_{B_1}\}$.
Let us proceed the same successive conjugate operations in Section 3 (see the outline of proof of Theorem 0.7) for $w_\tau=\tau_{j_1}\cdots\tau_{j_t}$. Then
after conjugations $\tau_{m_l}\tau_i\tau_{m_l}^{-1}$ of times at most $R$,
all $\tau_{m_l}$ in the word $w_\tau$ are moved to the last position:
$$\alpha_{j_1'}\tau_{j_1'}\cdots \alpha_{j_{t-n(m_l)}'}\tau_{j_{t-n(m_l)}'} \cdot \tau_{m_l}^{n(m_l)},$$
where each $|\alpha_{j_i}'|_{B_1}\le K$ and $n(m_l)$ is the total count of $\tau_{m_l}$ appears in $\tau_{j_1}\cdots\tau_{j_t}$. Moreover, by Lemma 3.2,
$$\alpha_{j_1'}\tau_{j_1'}\cdots \alpha_{j_{t-n(m_l)}'}\tau_{j_{t-n(m_l)}'} \cdot \tau_{m_l}^{n(m_l)}=\alpha_{m_l}'\tau_{j_1'}\tau_{j_2'}\cdots \tau_{j_{t-n(m_l)}'} \cdot \tau_{m_l}^{n(m_l)},$$
and $\alpha_{m_l}'$ has word length $\le P_0(K(R-1)+R)$.

Repeating the process for $\tau_{m_l-1},\dots,\tau_{m_{l-1}+1}$ one after another, we derive that
$$w_\tau=\alpha' \tau_{m_{l-1}+1}^{n(m_{l-1}+1)}\cdots\tau_{m_l}^{n(m_l)},\quad |\alpha'|_{B_1}\le R\cdot P_0(K(R-1)+R), \tag{4.5}$$
which implies (3.6) for $\Gamma_l$.

Now let us consider $1\to \Gamma_{l}\to \Gamma_{l-1}\to \Lambda_{l-1}/\Lambda_{l}\to 1$. By (4.5), Theorem 1.3 still holds for $\Gamma_l$ (not necessarily knowing $\Gamma_l$ is virtually nilpotent), where the ``graded'' basis formed by $B_1$ and $T_l$ of $\Gamma_l$ is associated to
$$\Gamma_l\vartriangleright N \vartriangleright\cdots\vartriangleright N_k \vartriangleright \{e\}.$$ On the other hand, for any $\bar\tau_i$ of infinite order, the eigenvalues of conjugation $\bar\rho(\bar \tau_j)\in \op{Aut}(\Gamma_l/N)$ on the free part have norm $1$. It is not difficult to see that the proof of Lemma 3.2 still works for $(\Gamma_{l-1},\Gamma_l)$.
By the proof of (4.5) again, (3.6) holds for $\Gamma_{l-1}$. By induction, the proof of Lemma 3.4 is complete.
\qed\enddemo

\vskip4mm

\head 5. Proof of Theorem 0.9
\endhead

\vskip4mm

We first recall the isoperimetric functions for a finitely presented group. Let $G$ be a group given by a finite presentation, $G=\left<S | \Cal R\right>$,
where $S=\{\gamma_1, \cdots, \gamma_m\}$ is a finite symmetric set generating $G$, and  $\Cal R\subset F_S$ is the set of relators, and $F_S$ is the free group generated by $S$. Because by definition $G$ is the quotient of $F_S$ by the normal closure of $\Cal R$ in $F_S$, a word $w=s_1\cdots s_n\in F_S$, where $s_i\in S$, is represented the identity $e\in G$ if and only if there are reduced words $v_i$ on $S$ such that $w$ can be written in $F_S$ as
$$w=\prod_{i=1}^k v_ir_i v_i^{-1}, \quad r_i \text{ or } r_i^{-1}\in \Cal R.\tag{5.1}$$

A function $f: \Bbb N \to \Bbb N$ is called an isoperimetric function of $\left<S | \Cal R\right>$ for any word of length $n$ that equals to the identity in $G$, if $w$ can be written in (5.1) with at most $k\le f(n)$ relators.
It is a analogue of the isoperimetric inequality in differential geometry, because when looking at the van Kampen diagram with boundary label $w$, then its boundary is of length $|w|_S$ and the numbers of relators $r_i$ in (5.1) corresponds to its ``area''.
The smallest isoperimetric function of a finite presentation is called the Dehn function, and it is well known that Dehn functions corresponding
to different finite presentations of the same group are equivalent (\cite{MO}).
A finitely presented group $G$ is said to satisfy a polynomial
isoperimetric inequality if there is a finite presentation whose isoperimetric function can be chosen to be a polynomial $Cn^d$, and thus so does any other finite presentation of $G$.

A function $f$ is called an isodiametric function of a finite presentation $\left<S | \Cal R\right>$ if for
every number $n$ and every word $w$ of length $n$ which is equal to $e$ in $G$,
$w$ is a product of conjugates in (5.1) such that the length of each $v_i$ is bounded by $f(n)$, or equivalently, there is
a van Kampen diagram with boundary label $w$ and diameter $\le f(n)-K$, where $K=\max\{|r|_{F_S}, \, r\in \Cal R\}$.
The following bound on the length of $r_i$ is used in the proof for (0.9.2).

\proclaim{Lemma 5.2}\rm{(\cite{Ge})}
If $f$ is an isoperimetric function for $\left<S | \Cal R\right>$, then $Kf(n)+n+K$ is an isodiametric function for $\left<S | \Cal R\right>$.
\endproclaim

The main technical result in proving Theorem 0.9 is an generalization of Theorem 0.7 that is applicable when $\Lambda$ is not virtually nilpotent. Let $1\to N\to \Gamma \overset{\pi}\to\longrightarrow \Lambda\to 1$ be a exact sequence as (0.6). Because $N$ is nilpotent, $h_w(N)=0$. Thus, the growth rate $h_w(\Gamma, S)$ of $\Gamma$ is determined by that of $\Lambda$, $h_w(\Lambda, \bar S)$, and
the distortion function of $N$ in $\Gamma$, $$\Delta_{N}^{\Gamma}(R)=\sup_{g\in N, |g|_{S}\le R} |g|_{N\cap S}.\tag 5.3$$
Based on the proof of Theorem 0.7, we will prove that for any ``horizontal'' lifting subgroup $\Gamma_0=\left<S_0\right>$ of $\Lambda=\left<\bar S_0\right>$, $h_w(\Gamma_0,S_0)$ coincides with that of the whole group, and then prove that the distortion function $\Delta_{N}^\Gamma$ is polynomial if the isoperimetric function of $\Lambda$ is polynomial.

\proclaim{Theorem 5.4} Let $\Gamma$ be a finitely generated group, and let $N$ be a normal nilpotent subgroup of $\Gamma$. Assume $(\Gamma,N)$ satisfies (0.6) and
the following two conditions:

\noindent \rm{(5.4.1)} For each $1\le h\le k$, $\op{Im}(\bar \rho_h)$ is virtually nilpotent.

\noindent \rm{(5.4.2)} There are $\gamma_1,...,\gamma_m\in \Gamma-N$
such that $\{\pi(\gamma_i)\}$ generates
$\Lambda$, and $\{\bar \rho_h(\pi(\gamma_i))\}$ contains a graded symmetric generating set for $\op{Im}(\rho_h)$
such that if $\pi(\gamma_i)$ has an infinite order, then eigenvalues of $\bar \rho_h(\pi(\gamma_i))$
have norm one.

Let $S$ be a finite generator system of $\Gamma$.

\noindent \rm{(5.4.3)}  If $S_0\subset S$ is a subset such that $\pi(S_0)=\pi(S)$,  then $h_w(\Gamma,S)=h_w(\left<S_0\right>,S_0)$.

\noindent \rm{(5.4.4)} If $\Lambda$ satisfies a polynomial isoperimetric inequality,
then $h_w(\Gamma,S)=h_w(\Lambda,\pi(S))$.
\endproclaim

Because a virtually nilpotent group admits a polynomial isoperimetric function, Theorem 5.4 is a natural extension of Theorem 0.7.

\demo{Proof of Theorem 5.4}

By Remark 4.4, Lemma 3.2 and (3.8) still hold for a graded generator system $B_1$ of $N$ and ${\gamma_1,\dots,\gamma_m}$ given in (5.4.2), where $\Lambda$ is not virtually nilpotent.
Let $S$ be given in Theorem 5.4. Since $h(\Gamma,S\cup B_1)\ge h(\Gamma, S)$, without loss of generality, we assume that $S$ contains $B_1$. Since the two generating sets
can be represented by each other, (3.8) also holds for a word $\gamma_{j_1}\cdots\gamma_{j_t}$ consisting of elements in $S$ and $\alpha_p\in S\setminus S_0$ or $\alpha_p\in B_1$.

(5.4.3) By Lemma 1.3 and Lemma 3.2 we derive
$$|\Gamma_0(R, S_0)|\le |\Gamma(R,S)|\le P(s)\cdot |\Gamma_0(R, S_0)|,$$
and thus (5.4.3) follows:
$$h_w(\Gamma,S)=h_w(\Gamma_0, S_0).$$

(5.4.4) Let $S_0\subset S$ be a subset whose projection $\bar S_0=\pi(S)$. Let us fix finite presentations, $\Lambda=\left<\bar S_0 | \Cal R\right>$, $\Gamma_0=\left<S_0 | \Cal R_0\right>$. For any word $w=\gamma_{i_1}\cdots\gamma_{i_R}\in \Gamma_0(R,S_0)\cap N$,  $\bar w=\bar\gamma_{i_1}\cdots\bar\gamma_{i_R}= 1\in \Lambda$.
Let $F_{\bar S_0}$ be the free group generated by $\bar S_0$, then there are reduced words $\bar v_i$ in $F_{\bar S_0}$ such that the word $\bar w\in F_{\bar S_0}$ can be written in the form as (5.1)
$$\bar w=\prod_1^n \bar v_i\bar r_i\bar v_i^{-1}, \qquad\bar r_i\text{ or }\bar r_i^{-1}\in \Cal R,$$
where $n$ is bounded by an isoperimetric function of $\Lambda$, and $|\bar v_i|_{\bar S_0}$ is bounded by an isodiametric function.

Accordingly, the word $w$ in the free group $F_{S_0}$ generated by $S_0$ can be written in
$$w=\prod_1^n v_i r_i v_i^{-1},$$
where $v_i$ and $r_i$ are the reduced words after replacing elements in $\bar S_0$ with the corresponding ones in $S_0$.
Then after projected to $\Gamma_0$, $r_i\in N$.
Let $B_1$ be the graded basis defined for $N$ in Section 3, and let
$$L=\max_{\bar r_i\in \Cal R}\{|r_i|_{B_1}\}, \qquad K=\max_{\bar r_i\in \Cal R}\{|\bar r_i|_{F_S}\}.$$
By the assumption in (5.4.4), $n$ is bounded by an isoperimetric function $f(R)$ which is polynomial in $R$. By Lemma 5.2, each $v_i$ has word length bounded by a polynomial $K(f(R)+1)+R$.
Then by (3.8),
$$|w|_{B_1}=|\prod_{i=1} v_ir_iv_i^{-1}|_{B_1}\leq LP_0(Kf(R)+R+K)\cdot f(R).\tag{5.5}$$
Because any word $\Gamma(R,S)$ can be rewritten by elements in $S_0\cup B_1$, by Lemma 3.2 and (5.5), the distortion function of $N$ in $\Gamma$
 is bounded by
$$\Delta_{N}^{\Gamma}(R)\le LP_0(Kf(R)+R+K)\cdot f(R)+P_0(R). \tag 5.6$$

Since $N$ is of polynomial growth, by (5.6) $|\Gamma(R,S)\cap N|$ is also bounded by a polynomial in $R$.
Because
$$\Lambda(R,\bar S_0)\le |\Gamma(R,S)|\le |\Gamma(2R,S)\cap N|\cdot |\Lambda(R, \bar S_0)|,$$ we derive $h_w(\Gamma, S)=h_w(\Lambda, \bar S)$.
\qed\enddemo

Let us consider a compact $n$-manifold $M$ satisfying
$$\op{Ric}_M\geq -(n-1), \quad \op{diam}(M)\leq d, \quad h(M)<\epsilon(n,d),$$
where $\epsilon(n,d)$ is given in Lemma 2.1.  Let $\Gamma=\pi_1(M)$, and $N,\Lambda$ are given in Theorem 1.12. As discussed in the proof of Theorem 0.3, there is a normal nilpotent subgroup $N$ of $\Gamma$, such that $h_w(\Lambda,\Lambda(2d+1))=0$, where $\Lambda=\Gamma/N$.

The groups $N,\Gamma,\Lambda$ satisfy the conditions in Theorem 5.4. Indeed, (5.4.2) follows directly from (2.1.2).  Because $h_w(\op{Im}(\bar \rho_{h}))\le h_w(\Lambda,\Lambda(2d+1))=0$ and $\op{Im}(\bar \rho_h)$ is a linear group, by Lemma 1.6 we conclude (5.4.1) that $\op{Im}(\bar \rho_{h})$ is virtually nilpotent for every $1\le h\le k$.

Now we are ready to prove Theorem 0.9.

\demo{Proof of Theorem 0.9}

Continue from the above discussion.

\noindent (0.9.1) If $\Gamma=N\rtimes_{\psi} \Lambda$, then there is an injective homomorphism, $\tau : \Lambda\to \Gamma$, such that $\tau(\Lambda)= \left<\tau(\Lambda)\right>$.  And thus $h_w(\tau(\Lambda),\tau(\Lambda(2d+1)))=h_w(\Lambda, \Lambda(2d+1))=0$. By Theorem 5.4 and Lemma 1.2,
$$ h(M) \leq  h_w(\Gamma,\Gamma(2d+1))\overset{\op{(5.4.3)}}\to =h_w(\tau(\Lambda),\tau(\Lambda)(2d+1))=0.$$

\noindent (0.9.2) If $\Lambda$ has a polynomial isoperimetric function, by Theorem 5.4 and (2.1.1),
$$ h(M)\leq h_w(\Gamma, \Gamma(2d+1))
\overset{\op{(5.4.4)}}\to=h_w(\Lambda, \Lambda(2d+1))=0.\qed$$
\enddemo

\vskip4mm

\head 6. Manifolds of Small Volume Entropy
\endhead

\vskip4mm

Let $M$ be a compact $n$-manifold satisfying
$$\op{Ric}_M\ge -(n-1),\quad \op{diam}(M)\le d.$$
By Bishop volume comparison, $0\le h(M)\le n-1$, and $h(M)=n-1$
iff $M$ is a hyperbolic manifold (\cite{LW}). In [CRX], we generalize
the maximal volume entropy rigidity to that if $h(M)$ is almost $n-1$
i.e., there is $\epsilon(n,d)>0$ such that if $h(M)\ge (n-1)-\epsilon(n,d)$, then
$M$ is diffeomorphic to a hyperbolic manifold by an $\Psi(\epsilon|n,d)$-GHA,
where $\Psi(\epsilon|n,d)\to 0$ as $\epsilon\to 0$ while $n$ are $d$ are fixed.

On the other hand, $h(M)=0$ implies that $\pi_1(M)$ is amendable
(\cite{Ga}, \cite{Ta}, cf. Corollary 3.6 \cite{DZ}) i.e. $\pi_1(M)$ admits an invariant additive probability measure, but the converse does not hold
(see Lemma 1.7). Note that finitely presented groups are not necessarily amendable (\cite{OS}).

In view of the above, a natural problem is whether $\pi_1(M)$ is amendable when $h(M)$ is small; while the gap vanishing phenomena suggested by the
Nilpotency Conjecture 0.2 is still open. We have the following positive answer:

\proclaim{Lemma 6.1} Given $n, d$, there exists $\epsilon(n, d)>0$ such that
if a compact $n$-manifold $M$ satisfies
$$\op{Ric}_M\ge -(n-1), \quad d\ge \op{diam}(M),\quad \epsilon(n,d)>h(M),$$
then $\pi_1(M)$ is amendable.
\endproclaim

\demo{Proof} Let $\epsilon(n, d)>0$ be as in Lemma 2.1. Then by (2.1.1), $h_w(\Lambda, \Lambda(2d+1))=0$ which implies $\Lambda$ is amenable.
And by (0.6), there is $1\to N\to \Gamma\to \Lambda\to 1$, where $N$ is a normal nilpotent subgroup. Since a normal extension of an amenable group is amenable (cf. \cite{Pi}), we have that $\Gamma$ is amenable.
\qed\enddemo

We point it out that in the proof of (0.8.2), we present a direct proof for $\Gamma$
to not contain any free group of two generators. This fact also follows from
Lemma 6.1 and a known fact that any amenable group contains no $2$-generated free groups
(cf. \cite{Wa}).

Next, we will verify the gap vanishing volume entropy property for
manifolds satisfying additional conditions which may not be enough for one
to verify the Nilpotency Conjecture 0.2.

\proclaim{Lemma 6.2}
Given $n, d, v>0$ there exists $\epsilon(n,d,v)>0$ such that
if a compact $n$-manifold satisfies
$$\op{Ric}_M\ge -(n-1),\quad d\ge \op{diam}(M),\quad \op{vol}(M)\ge v,\quad \epsilon(n,d,v)>h(M),$$
then $h(M)=0$.
\endproclaim

\demo{Proof}
Arguing by contradiction, assume a contradicting sequence, $M_i$, that satisfies $\op{Ric}_{M_i}\ge -(n-1)$, $\op{diam}(M_i)\le d$ and $0<h(M_i)\to 0$. Passing to a subsequence, we may assume a map,
$\phi_i: \Gamma_i\to G$, such that (1.8) holds. By Theorem 1.9, there is $\epsilon>0$
such that $\Gamma_i$ has a  normal subgroup $\Gamma_{i,\epsilon}$ such that
$K_i=\Gamma_i/\Gamma_{i,\epsilon}\cong K=G/G_0$. Because $\op{vol}(M_i)\ge v>0$,
$G_0=e$ i.e., $G$ is discrete. Consequently,  $\Gamma_{i,\epsilon}$ is at most finite
(actually trivial when $i$ large, see Lemma 1.10). Then $\Gamma_i(2d+1)$ is a finite symmetric
generating set (Lemma 1.2) and $h_w(G,G(2d+1))=h_w(K_i,K_i(2d+1))=
h_w(\Gamma_i,\Gamma_i(2d+1))$.  Again by Lemma 1.2,
$$h(M_i)\le h_w(G, G(2d+1))=\lim_{i\to \infty}h_w(\Gamma_i,\Gamma_i(2d+1))\le (2d+1)\lim_{i\to \infty}h(M_i)=0,$$
a contradiction.
\qed\enddemo

As seen in the proof of Theorem 0.9, if $\Delta_N^\Gamma(R)$ in (5.3) is bounded by a polynomial in $R$, then
$h(M)<\epsilon(n,d)$ implies that $h(M)=0$. Based on the work of \cite{CS}, \cite{Le}, the
desired property holds in the following situation.

\proclaim{Theorem 6.3} Given $d>0$, there is a constant $\epsilon(n,d)>0$ such that for a compact $n$-manifold $M$ satisfying $\op{Ric}_M\ge -(n-1)$ and $d\ge \op{diam}(M)$.

\noindent (6.3.1) If $M$ is a manifold without conjugate points, then
$h(M)=0$ or $h(M)\ge \epsilon(n,d)$;

\noindent (6.3.2) If $M$ is of non-positive sectional curvature,
then $M$ is flat or  $h(M)\ge \epsilon(n,d)$.
\endproclaim

\demo{Proof} Let $\epsilon(n, d)>0$ be as in Lemma 2.1.

\noindent (6.3.1) From \cite{CS} and \cite{Le} (also an unpublished work of Kleiner),
a compact manifold without conjugate point satisfies that
every solvable subgroup $W$ of $\Gamma=\pi_1(M)$ is a Bieberbach group and straight in $\Gamma$, i.e., the distortion $\Delta_W^\Gamma(R)$ is linear in $R$ with respect to any generator system $S$ of $\Gamma$ such that $S\cap W$ generates $W$. In particular, $N$ is in (0.6) satisfies $$\Delta_N^\Gamma(R)\le c(S)R,$$ where $c(S)$ is a constant that depends on $S$. As in the proof of Theorem 0.9 we conclude that $h(M)=0$ or $h(M)\ge \epsilon(n,d)$.

\noindent (6.3.2) By Lemma 6.1, $h(M)\ge \epsilon(n,d)>0$ or $\pi_1(M)$ is amenable. By \cite{An} and \cite{Av}, if $\pi_1(M)$ is amendable, then $\pi_1(M)$ is Biebebach, and thus $M$ is flat (because $\op{sec}_M\le 0$).
\qed\enddemo

\remark{Remark \rm 6.4} Let $M$ be a compact manifold without conjugate points. By \cite{FM}, $h(M)$
equals to the topological entropy of geodesic flows and they conjectured that if $M$ has a vanishing topological entropy, then  $M$ must be flat.
\endremark

\vskip4mm

\head 7. Index Bound In The Nilpotency Conjecture 0.2
\endhead

\vskip4mm

As discussed in Remark 0.4, the following problem is open:

\example{Problem 7.1} Does the following strengthened Nilpotency Conjecture 0.2 holds:
given $n, d>0$, are there constants, $\epsilon(n,d), w(n,d)>0$, such that if
$$\op{Ric}_M\ge -(n-1),\quad d\ge \op{diam}(M),\quad h(M)<\epsilon(n,d),$$
then $\pi_1(M)$ is $w(n,d)$-nilpotent.
\endexample

From \cite{KW}, it can be seen that Problem 7.1 has a positive answer if either of the following conditions are satisfied:

\noindent (7.2.1) The diameter is sufficient small;

\noindent (7.2.1) $\pi_1(M)$ is finite.

In this section we prove that a positive answer of Milnor's problem implies that the assertion in Problem 7.1 holds for torsion free groups.

\proclaim{Theorem 7.3} If Milnor Problem 0.1 has a positive answer, then given $n, d>0$, there are constants, $\epsilon(n,d), w(n,d)>0$, such that if
$$\op{Ric}_M\ge -(n-1),\quad d\ge \op{diam}(M),\quad h(M)<\epsilon(n,d),$$
and $\pi_1(M)$ is torsion free,
then $\pi_1(M)$ is $w(n,d)$-nilpotent.
\endproclaim

Theorem 7.3 replies on the fact that if $N$ is a simply connected nilpotent Lie group, then
any discrete subgroup $\Gamma$ of  $\op{Aut}(N)$ such that $N/\Gamma$ is compact
satisfies $[\Gamma:\Gamma\cap N]\le w(n)$. In general, given any $C>0$, there is a virtually nilpotent group that is not $C$-nilpotent.

\example{Example 7.4} For any prime  $p$, let $\Gamma$ satisfy
$$0\to \Bbb Z_p \to \Gamma=\Bbb Z_p\rtimes_\varphi \Bbb Z \to \Bbb Z\to 0.$$
Then the generator $\alpha$ (resp. $\bar\gamma$) of $\Bbb Z_p$  (resp. of $\Bbb Z$) satisfies $\varphi( \gamma)(\alpha)=\gamma\alpha\gamma^{-1}=\alpha^2$, where $\gamma$ is a representative of
$\bar \gamma\in \Gamma$. The nilpotent subgroup is $\left<\gamma\right>$ whose index in $\Gamma$ is $p$. It is suspected in \cite{KW} that such groups cannot be a fundamental group of a manifold in Problem 7.1.
\endexample

\demo{Proof of Theorem 7.3}

Let $M$ be a compact $n$-manifold satisfying
$$\op{Ric}_M\geq -(n-1), \quad \op{diam}(M)\leq d, \quad h(M)<\epsilon(n,d),$$
where $\epsilon(n,d)$ is given in Lemma 2.1. Let $1\to N\to \Gamma\to \Lambda\to 1$ be the exact sequence (0.6) where
$\Gamma=\pi_1(M)$, and $N,\Lambda$ are given in Theorem 1.12.
As discussed in the proof of Theorem 0.3, $h_w(\Lambda)=0$ and the conjugate matrix $\rho_h(\bar\gamma)$ on $N_h/N_{h+1}$ has eigenvalues of norm $1$ if $\bar \gamma$ is of infinity order. If Milnor Problem 0.1 has a positive answer, then $\Lambda$ is virtually nilpotent.
Together with the finiteness of $\Lambda$ and the conjugate action $\bar\rho_h:\Lambda\to\op{Aut}(N_h/N_{h+1})$,
Theorem 7.3 will follow from the following Proposition 7.5 below.
\qed \enddemo

\proclaim{Proposition 7.5} Let $N,\Gamma,\Lambda$ be as in (0.6) and $N$ is torsion free. If

\noindent \rm{(7.5.1)} $\Lambda$ is virtually nilpotent, and

\noindent \rm{(7.5.2)} for each $\bar\gamma\in \Lambda$ with infinite order, all eigenvalues of $\bar\rho_h(\bar\gamma)$ has norm $1$,

\noindent then $\Gamma$ is $\omega(\Lambda, \bar \rho_h)$-nilpotent.
\endproclaim

The proof of Proposition 7.5 is based on an argument in \cite{Wo} (Proposition 4.4, Page 434), and Mal'cev's theorem on nilpotent groups.
Recall that a lattice in a Lie group G is a discrete subgroup $\Gamma$ with the quotient $G/\Gamma$ compact.

\proclaim{Theorem 7.6 (\cite{Ma})}
Any finitely generated torsion free nilpotent group can be embedded as
a lattice in a simply connected nilpotent Lie group $G$.
\endproclaim

\demo{Proof of Proposition 7.5}

Let $\Gamma_0$ be a subgroup of $\Gamma$ satisfying
$1\to N\to \Gamma_0\to \Lambda_0\to 1,$
where $\Lambda_0$ is a nilpotent subgroup of $\Lambda$ with index $c(\Lambda)<\infty$.

By Theorem 7.6, there is a connected simply connected nilpotent Lie group $D$ containing $N$ as a discrete subgroup with $D/N$ is compact.
For $\gamma\in \Gamma_0$, $\alpha\in N$, let $\rho(\gamma)(\alpha)=\gamma\alpha\gamma^{-1}$ and let $\rho_{\ast}(\gamma)$ be the induced automorphism of the Lie algebra $\Cal D$ of $D$. Then it is easy to see that after fixing a graded basis of $N$, the block matrix of $\rho_{*}(\gamma)$ around its diagonal coincides with the corresponding matrix of $\rho_h(\gamma)$, and thus the eigenvalues of $\rho_*(\gamma)$ consists of $\rho_h(\gamma)$'s. Moreover, by its Jordan form each matrix $\rho_*(\gamma)$ can be written in a product $A_\gamma\cdot B_\gamma$ of two matrices, where $A_\gamma$ is diagonal formed by all eigenvalues of $\rho_*(\gamma)$ whose norm equals $1$. Then $\{A_\gamma, \gamma\in \Lambda_0\}$ forms an abelian group in $\op{Aut}(\Cal D)$, whose closure $\bar A$ is compact. Since every element in $\bar A$ induces an automorphism of $D$ which preserves the lattice $N$, any connected subgroup of $\bar A$ is trivial. Therefore, we derive that $\bar A$ is not only compact but also discrete. Thus $\bar A$ is finite, and its order is bounded by a constant depending on $\{\rho_h,\, 1\le h\le k\}$.

Now let us consider the homomorphism $\Gamma_0\to A$, $\gamma\mapsto A_\gamma$. We claim that its kernel,
$$U=\{\gamma\in \Gamma_0,\, \rho_{\ast}(\gamma) \text{ has eigenvalue }1\},$$ is nilpotent.

By the claim, $$[\Gamma:U]\le[\Gamma:\Gamma_0]\cdot [\Gamma_0:U]=c(\Lambda)\cdot |\bar A|= c_1(\Lambda,\rho_h).$$ The proof of Proposition 7.5 is complete.

The claim follows from the same arguments as the proof of Proposition 4.4 in \cite{Wo}. Because in their case $\Lambda$ was assumed to be free abelian, we present a proof for completeness. Let us prove by induction on the dimension $m=\op{dim}(D)$.
If $m=1$, then for any $\gamma\in U$, $\alpha\in N$, $\gamma\alpha\gamma^{-1}=\alpha$, i.e., $N$ is in the center of $U$. Let us take a graded symmetric generating set of $\Lambda_0$, $\{\bar\gamma_1,\cdots, \bar\gamma_{n_1}, \cdots, \bar\gamma_{n_s}\}$, $\Lambda_{0,h}=\left<\bar\gamma_{n_{h-1}+1}\cdots, \bar\gamma_{n_s}\right>$ and $[\Lambda_0, \Lambda_{0,h}]\subset \Lambda_{0,h+1}$. Let $U_h=\left<\gamma_{n_{h-1}+1},\cdots, \gamma_{n_s}, \alpha_0\right>\cap U$, where $\alpha_0$ is a generator of $N$.
Then $[U, U_h]\subset U_{h+1}$, and $U_{s+1}=1$. Thus $U$ is nilpotent.

Assume the claim for $m=k$. For $m=k+1$,  take a  central series of $D$,
$$D=D_1\vartriangleright D_2\vartriangleright\cdots\vartriangleright D_{k+1}\vartriangleright 1,$$
such that $[D, D_h]\subset D_{h+1}$ and $D_h/D_{h+1}$ is of dimension $1$.
Since for any $\gamma\in U$, $\rho_{\ast}(\gamma)$'s every eigenvalue equals to $1$, there is a basis of $\Cal D$ such that the matrix $\rho_*(\gamma)$ is transformed to be upper triangular for every $\gamma\in U$ at the same time. Thus the central series can be chosen such that $\rho(U)$ preserves each $D_h/D_{h+1}$ and $D_h/(N\cap D_{h+1})$ is compact.
In particular, $N\cap D_{k+1}$ is cyclic and lies in the center of $U$. By the inductive assumption for $(U/(N\cap D_{k+1}), N/(N\cap D_{k+1}))$, where $\op{dim}(D/D_{k+1})=k$, $U/(N\cap D_{k+1})$ is nilpotent.
Applying the same arguments above for the case $\dim(D)=1$ to the short exact sequence,
$$1\to N\cap D_{k+1}\to U\to U/(N\cap D_{k+1})\to 1,$$
we conclude that $U$ is nilpotent. \qed\enddemo


\vskip4mm



\Refs
\nofrills{References}
\widestnumber\key{APS1}

\vskip3mm


\ref
\key An
\by M. T. Anderson
\pages 269-278
\paper On the fundamental group of nonpositively curved manifolds
\jour Math. Ann
\vol 276
\yr 1987
\endref

\ref
\key Av
\by A. Avez
\pages 188-191
\paper Vari\'et\'es riemanniennes sans points focaux
\jour C. R. Acad. Sci. Paris
\vol I 270
\yr 1970
\endref

\ref
\key BH
\by M. Bucher; P. de la Harpe
\pages 811-815
\paper Free products with amalgamation and HNN-extensions which are of uniformly exponential growth
\jour Math. Notes
\vol 67
\issue 6
\yr 2000
\endref



\ref
\key Ce
\by L. Di Cerbo
\pages 193-199
\paper A gap property for the growth of closed 3-manifold groups
\jour Geometriae Dedicata
\vol 143
\yr 2009
\endref

\ref
\key CG
\by J. Cheeger; D. Gromoll
\pages 119-128
\paper The splitting theorem for manifolds of nonnegative Ricci curvature
\jour J. Diff. Geom.
\vol 6
\yr 1971
\endref


\ref
\key CRX
\by L. Chen, X. Rong; S. Xu
\pages
\paper Quantitative Volume Space Form Rigidity Under Lower Ricci Curvature Bound
\jour J. Diff. Geom.
\vol
\yr To appear
\endref

\ref
\key Ch
\by C. Chou
\pages 396-407
\paper Elementary amenable groups
\jour Illinois J. Math.
\vol 24
\yr 1980
\endref

\ref
\key CN
\by T. Colding; A. Naber
\pages 1173-1229
\paper Sharp H\"older continuity of tangent cones for spaces with a lower Ricci curvature bound and applications
\jour Ann. of Math.
\vol 176 (2)
\yr 2012
\endref

\ref
\key CZ
\by D. J. Collins; H. Zieschang
\paper Combinatorial group theory and fundamental groups
\inbook Algebra, VII
\bookinfo Encyclopaedia Math. Sci.
\publ Springer
\publaddr Berlin
\vol 58
\pages 1-166
\yr 1993
\endref


\ref
\key CS
\by C. Croke; V. Schroeder
\pages 161-175
\paper The fundamental group of compact manifolds without conjugate points
\jour Comment. Math. Helvetici
\vol 61
\yr 1986
\endref



\ref
\key DZ
\by  M. Do Carmo; D. Zhou
\pages 1391-1401
\paper Eigenvalue estimate on complete noncompact Riemannian manifolds and applications
\jour Trans. Amer. Math. Soc.
\vol 351
\yr 1999
\endref


\ref
\key FM
\by A. Freire;  R. Ma\~n\'e
\pages 375-392
\paper On the entropy of the geodesic flow in manifolds without
conjugate points
\jour  Invent. Math.
\vol 69
\yr 1982
\endref

\ref
\key FY
\by K. Fukaya; T. Yamaguchi
\pages 253-333
\paper The fundamental groups of almost non-negatively curved manifolds
\jour Ann. of Math.
\yr 1992
\vol 136
\endref

\ref
\key Ga
\by M.E. Gage
\pages897-912
\paper Upper bounds for the first eigenvalue of the Laplace-Beltrami operator
\jour Indiana Univ. Math. J
\vol 29
\yr 1981
\endref

\ref
\key Ge
\by Steve M. Gersten
\paper Isoperimetric and isodiametric functions of finite presentations
\inbook Geometric group theory. Volume 1.
\ed Niblo, Graham A. et al.
\procinfo Proceedings of the symposium held at the Sussex University, Brighton (UK), July 14-19, 1991
\publ Cambridge: Cambridge University Press
\bookinfo Lond. Math. Soc. Lect. Note Ser. 181
\yr 1993
\pages 79-96
\endref

\ref
\key Gri1
\by R. Grigorchuk
\pages 31-33
\paper On Milnor's problem of group growth
\jour Dikl. Ak. Nauk SSSR
\vol 271
\yr 1983
\endref

\ref
\key Gri2
\by R. Grigorchuk
\pages 939-985
\paper The growth degrees of finitely generated groups and the theory of invariant means
\jour Izv. AKad. Nauk SSSR. Ser. Math.
\vol 48
\yr 1984
\endref

\ref
\key Gri3
\by R. Grigorchuk
\pages 194-214
\paper On the degrees of p-groups and torsion-free groups
\jour Math. Sbornik
\vol 126
\yr 1985
\endref

\ref
\key Gri4
\by R. Grigorchuk
\pages
\paper Milnor's problem on the growth of groups and its consequences
\jour arXiv: 1111.0512
\vol
\yr 2013
\endref


\ref
\key GP
\by R. Grigorchuk; I, Pak
\pages 251-272
\paper  Groups of intermediate growth: an introduction
\jour Enseign. Math. (2)
\vol 54
\issue 3-4
\yr 2008
\endref

\ref
\key Gr
\by M. Gromov
\pages  53-73
\paper Groups of polynomial growth and expanding maps
\jour Inst. Hautes, Etudes Sci. Publo Math
\vol 53
\yr 1981
\endref


\ref
\key KW
\by V. Kapovitch; B. Wilking
\pages
\paper Structure of fundamental groups of manifolds with Ricci curvature bounded below
\jour Preprint
\vol
\yr 2011
\endref

\ref
\key Le
\by N. Lebedeva
\paper On the fundamental group of a compact space without conjugate points
\jour PDMI Preprint
\yr 2002
\endref


\ref
\key LW
\by F. Ledrappier; X. Wang
\pages 461-477
\paper An integral formula for the volume entropy with application to rigidity
\jour J. Diff. Geom.
\vol 85
\yr 2010
\endref




\ref
\key Ma
\by A. I. Mal'cev
\pages 9-32
\paper On a class of homogeneous spaces
\jour Izvestiya Akad. Nauk. SSSR. Ser. Mat.
\vol 13
\yr 1949
\endref

\ref
\key Man
\by A. Manning
\pages 567-573
\paper Topological entropy for geodesic flows
\jour Ann. of Math
\vol 2
\yr 1979
\endref

\ref
\key MO
\by K. Madlener; F. Otto
\paper Pseudo-natural algorithms for the word problem for finitely presented monoids and groups
\jour J. Symbolic Computation
\vol 1
\yr 1989
\pages 383-418
\endref


\ref
\key Mi1
\by J. Milnor
\pages 1-7
\paper A note on curvature and fundamental group
\jour J. Diff. Geom.
\vol 2
\yr 1968
\endref

\ref
\key Mi2
\by J. Milnor
\pages 447-449
\paper Growth of finitely generated solvable groups
\jour J. Diff. Geom.
\vol 2
\yr 1968
\endref

\ref
\key Mi3
\by J. Milnor
\pages 685-686
\paper Problem 5603
\jour Amer. Math. Monthly
\vol 75
\yr 1968
\endref

	

\ref
\key OS
\by A. Olshanskiim; M. Sapir
\pages 43-169
\paper Non-amenable finitely presented torsion-by-cyclic groups
\jour Publ. Math. Inst. Hautes \'Etudes Sci.
\vol 96
\yr 2002
\endref

\ref
\key Os1
\by D. Osin
\pages 907-918
\paper The entropy of solvable groups
\jour Ergod. Th. \& Dynam. Sys.
\vol 23
\yr 2003
\endref

\ref
\key Os2
\by D. Osin
\pages 133-151
\paper Algebraic Entropy of Elementary Amenable Groups
\jour Geometriae Dedicata
\vol 107
\yr 2004
\endref

\ref
\key PR
\by J. Pan; X. Rong
\pages
\paper Ricci curvature and isometric action with scaling nonvanishing property
\jour In preparation
\vol
\yr
\endref


\ref
\key Pe
\by P. Petersen
\pages
\paper Riemannian Geometry (second edition)
\jour Springer-Verlag, New York
\vol
\yr 2006
\endref

\ref
\key Pi
\by J. Pier
\pages
\paper Amenable locally compact groups
\jour New York: Wiley
\vol
\yr 1984
\endref




\ref
\key Ro
\by X. Rong
\pages 193-298
\paper Convergence and collapsing theorems in
Riemannian geometry
\inbook Handbook of Geometric Analysis II
\bookinfo ALM 13
\publ Higher Education Press and International Press
\publaddr Beijing-Boston
\yr 2011
\endref


\ref
\key Ta
\by M. Taylor
\pages 773-793
\paper Lp-estimates on functions of the Laplace operator,
\jour Duke Math. J.
\vol 58
\yr 1989
\endref

\ref
\key Ti
\by J. Tits
\pages 250-270
\paper Free subgroups in linear groups
\jour Journal of Algebra
\vol 20
\yr 1972
\endref

\ref
\key Wa
\by S. Wagon
\book The Banach-Tarski paradox
\bookinfo reprint of the 1985 original
\publ Cambridge Univ. Press
\publaddr Cambridge
\yr 1993
\endref

\ref
\key Wo
\by J. Wolf
\pages 421-446
\paper Growth of finitely generated solvable groups and curvature of Rieamannian manifolds
\jour J. Diff. Geom.
\vol 2
\yr 1968
\endref

\endRefs
\enddocument